\documentclass[12pt,reqno,a4paper]{amsart}
\usepackage{verbatim,rcs,cite}
\usepackage{amssymb,amsfonts,latexsym,mathrsfs}
\usepackage{amsmath}
\usepackage{amsthm}
\usepackage{mathrsfs}
\usepackage[all]{xy}
\usepackage{pstricks}
\usepackage{hyperref}

\newcommand{\tmop}[1]{\ensuremath{\operatorname{#1}}}
\newcommand{\assign}{:=}
\newcommand{\nin}{\not\in}

\newtheorem{theorem}{Theorem}[section]
\newtheorem{lemma}[theorem]{Lemma}
\newtheorem{proposition}[theorem]{Proposition}
\newtheorem{corollary}[theorem]{Corollary}

\theoremstyle{remark}

\newtheorem{remark}[theorem]{Remark}
\theoremstyle{definition}
\newtheorem{definition}[theorem]{Definition}

\title{Topology of eigenspace posets for unitary reflection groups}

\date{}

\author{Justin Koonin}

\dedicatory{\upshape
School of Mathematics and Statistics \\
University of Sydney, NSW 2006, Australia\\[.5em]
{\it E-mail address:} \ \texttt{justin.koonin@sydney.edu.au} }

\begin{document}




\thanks{\noindent{The author is supported by ARC Grant \#DP110103451 at the University of Sydney.}}

\begin{abstract}
The eigenspace theory of unitary reflection groups, initiated by Springer
  and Lehrer, suggests that the following object is worthy of study: the poset of eigenspaces of elements of a unitary reflection group, for a fixed eigenvalue, ordered by the reverse of inclusion. We investigate topological properties of this poset. The new results
  extend the well-known work of Orlik and Solomon on the lattice of intersections of hyperplanes.
\end{abstract}

\maketitle

\section{Introduction}

Let $V$ be a complex vector space of finite dimension, and $G \subseteq GL(V)$ a unitary reflection group in $V$.  Denote by $\mathcal{A}(G)$ the set of
reflecting hyperplanes of all reflections in $G$, and
$\mathcal{M}_{\mathcal{A}(G)}$ the hyperplane complement -- that is, the smooth
manifold which remains when all the reflecting hyperplanes are removed from
$V$. There is an extensive literature studying the topology of
$\mathcal{M}_{\mathcal{A}(G)}$ (\cite{Arnold1969}, \cite{Brieskorn1973}, \cite{OrSo1980},
\cite{OrSo1980-2}, \cite{Lehrer1995}, \cite{BlLe2001}).

In particular, Orlik and Solomon \cite[Corollary 5.7]{OrSo1980} showed that $H^{\ast} (\mathcal{M}_{\mathcal{A}(G)}, \mathbb{C})$ is determined (as a graded representation of $G$) by the poset $\mathcal{L}(\mathcal{A}(G))$ of intersections of the hyperplanes in $\mathcal{A}(G)$.

It is well-known that the poset $\mathcal{L}(\mathcal{A}(G))$ is a geometric
lattice. Hence it is Cohen-Macaulay, and its reduced homology
vanishes except in top dimension. The poset $\mathcal{L}(\mathcal{A}(G))$ is
also known to coincide with the poset of fixed point subspaces (or
1-eigenspaces) of elements of $G$ (see Theorem \ref{proposition:fixedpoint}). Springer and Lehrer (\cite{Springer1974}, \cite{LeSp1999}, \cite{LeSp1999-2}) developed a general theory of eigenspaces for unitary reflection groups.
The purpose of this paper is to study topological properties of
generalisations of $\mathcal{L}(\mathcal{A}(G))$ for arbitrary eigenvalues.

Namely, let $\zeta$ be a complex root of unity, and $g$ be an element of $G$. Define $V (g, \zeta)
\subseteq V$ to be the $\zeta$-eigenspace of $g.$ That is, $V (g, \zeta) := \{v \in V \mid g
v = \zeta v\}.$  Let $\mathcal{S}_{\zeta}^V (G)$ be
the set $\{V (g, \zeta) \mid g \in G\}$, partially ordered by the reverse of
inclusion. More generally, if $\gamma \in N_{GL(V)} (G)$ (the
normaliser of $G$ in $GL(V)$) and $\gamma G$ is a reflection
coset, we may define $\mathcal{S}_{\zeta}^V (\gamma G)$ to be the set $\{V (x, \zeta)
\mid x \in \gamma G\}$, partially ordered by the reverse of inclusion.  This is a linear analogue of the poset of
$p$-subgroups of a group $G$ first studied by Quillen \cite{Quillen1978}.  The study of this poset was first suggested by Lehrer (see \cite[Appendix C, p.270]{LeTa2009}).

The main result of this paper is the following:

\begin{theorem}\label{theorem:CMall}
  Let $G$ be a unitary reflection group acting on $V$. Then $\mathcal{S}_{\zeta}^V (G)$ is Cohen-Macaulay over
  $\mathbb{Z}$.
\end{theorem}

The proof of this result depends on the following theorem, whose proof appears in 
\cite{Koonin2012-2} and \cite{Koonin2012}. The imprimitive reflection groups $G (r, p, n)$
are defined in Section \ref{subsection:classification}. 

\begin{theorem}\label{theorem:CMimprimitive}
  Let G be any imprimitive unitary reflection group $G
  (r, p, n)$, and $\gamma$ any element of $N_{GL (V)} (G)$ . The poset $\mathcal{S}_{\zeta}^V (\gamma G)$ is
  Cohen-Macaulay over $\mathbb{Z}$.
\end{theorem}

The proof of Theorem \ref{theorem:CMimprimitive} is combinatiorial in flavour, relying on the establishment of an isomorphism between the posets $\mathcal{S}_{\zeta}^V (\gamma G)$ and certain subposets of Dowling lattices.  The ideas used are very different from the invariant theoretic methods used in this paper, and for this reason the two papers have been separated.

The strategy for proving Theorem \ref{theorem:CMall} is to reduce the
statement to a computation of the reduced homology of a small number of
posets. All except two of these computations are somewhat routine and were performed
using MAGMA (see \cite{BoWiCa1997}, \cite{MAGMA}) and the GAP package
    Simplicial Homology (see \cite{GAP42008}, \cite{GAP4} for GAP and \cite{DuHeSaWe2003}, \cite{SimplicialHomology} for the Simplicial Homology package). The other two posets, associated with the
complexification of the real reflection group $E_8$, are very large. Additional techniques in computational algebraic topology were needed to
compute their reduced homology. A second paper \cite{KoJu2012} describing
these technques and their applications, written jointly with Mateusz Juda from the Jagiellonian University in Krakow, will appear shortly.

A third forthcoming paper \cite{Koonin2012-4} will consider an application of the
main result in this paper to the representation theory of unitary reflection
groups.

\section{Preliminaries}

\subsection{Notation}

Throughout this paper, $G$ will denote a group, and $V$ a complex vector space
of dimension $n$. If $H$ is a subgroup of $G$, write $N_G (H)$ for the
normaliser of $H$ in $G$. If $G \subseteq GL(V)$ and $x \in N_{GL(V)} (G)$, denote the
centraliser of $x$ in $G$ by $C_G (x)$. If $g \in GL(V),$ write
$\tmop{Fix} g =\{v \in V \mid gv = v\}$ for the fixed-point subspace of $g$. 
If $A$ and $B$ are two sets, the set of elements in
$A$ but not $B$ is denoted $A \backslash B$. If $x \in \tmop{End} (V)$, let $V (x, \zeta)$ be the $\zeta$-eigenspace of
$x$ acting on $V$.

\subsection{Partially ordered sets}

Let $P$ denote a poset. In this paper, all posets will be assumed to be
finite. If $a, b \in P$, {\em b covers a} if $a < c
\leqslant b$ implies $c = b$. A {\em lower order ideal} of $P$ is a
subset $I \subseteq P$ such that if $x \in I$ and $y \leqslant x$ then $y \in
I$. An {\em upper order ideal} of $P$ is a subset $I$ such that if $x
\in I$ and $y \geqslant x$ then $y \in I.$ \ If $x \in P$, denote by
$P_{\leqslant x}$ the subposet of $P$ consisting of elements less than or
equal to $x$. Similarly define $P_{< x}, P_{\geqslant x}, P_{> x}$. Define
the closed interval $[x, y] \assign \{z \in P \mid x \leqslant z \leqslant
y\}$. Similarly, define the open interval $(x, y) \assign \{z \in P \mid x <
z < y\}.$ If $P$ has a unique minimal element, denote this element
$\hat{0}$. An {\em atom} $x \in P$ is an element which covers
$\hat{0}$. Similarly, if $P$ has a unique maximal element, denote this
element $\hat{1}$. A {\em coatom} $x \in P$ is an element which is
covered by $\hat{1}$.

A {\em chain} is a
  poset in which any two elements are comparable. Define the
 {\em length} of a chain $C$ to be $\left| C \right| - 1,$ and write this
  number as $l (C) .$ The length of a poset $P$ is defined to be:
\[
l(P)\assign \max \{l (C) \mid C \mbox{ a chain in $P$}\}.
\]

A poset $P$ is
 {\em ranked} (or {\em pure}) of rank $n$ if all maximal
  chains in $P$ have the same length $n$. The {\em rank function} $r : P \longrightarrow \mathbb{N}$ of a ranked poset $P$ is defined inductively by $r(x) = 0$ if $x$ is a minimal element of $P$, and $r(y) = r(x) + 1$ if $y$ covers $x$.

Denote the {\em product} of two posets $P$ and $Q$ by $P \times Q$, and the
{\em join} by $P \ast Q$ (see \cite[\S3.2]{Stanley1997} for definitions).

The poset $P$ is a {\em join semilattice} if $x \vee y$ exists for all
$x, y \in P$ (which implies that $P$ has a unique maximum element $\hat{1}$). The poset $P$ is a {\em meet semilattice} if $x \wedge y$
exists for all $x, y \in P$ (which implies that $P$ has a unique minimum element $\hat{0}$). If $P$ is both a join semilattice and a meet semilattice, then $P$ is called a {\em lattice}.

A lattice $L$ is said to be {\em atomic} if every element of $L$ is the
join of atoms (by convention, the unique minimal element $\hat{0}$ is the join of the
empty set of atoms). A lattice $L$ is
{\em semimodular} if it is ranked, and if the associated rank function
$r$ of $L$ has the property that $r(x \vee y) + r (x \wedge y) \leqslant r (x) + r (y)$ for all
  $x,y \in L$. A lattice is {\em geometric} if it is both atomic and semimodular.

\subsection{Poset topology}

For the basics of poset topology, see \cite{Wachs2004}.

The {\em order complex} of a poset $P$, denoted $\Delta (P)$, is the
abstract simplicial complex whose vertices are the elements of $P$ and whose
simplices are the chains of $P$.

The simplicial homology of $P$ over a commutative ring $\mathbb{A}$, denoted
$H (P, \mathbb{A})$, is the simplical homology of the order complex $\Delta
(P)$.  Where the ring $\mathbb{A}$ is clear, we may denote this simply by $H (P)$. Similarly, reduced homology can be defined for $P$, and is denoted
$\tilde{H} (P, \mathbb{A}).$ The {\em n}th homology and reduced
homology modules are denoted $H_n (P, \mathbb{A})$ and $\tilde{H}_n (P,
\mathbb{A})$ respectively.

A poset is said to be {\em connected} (resp. {\em contractible}) if the associated simplicial
complex is connected (resp. contractible).

\subsection{Cohen-Macaulay posets}

The property of being Cohen-Macaulay has its origins in commutative algebra,
in the context of Cohen-Macaulay rings. For more details, see \cite[Lecture 4]{Wachs2004}.

Let $C = (x_0 < \ldots < x_i)$ be a chain
  of $P$. Define the {\em link} of $C$ in $P$ by
\[
\tmop{lk}_P (C) =\{z \in P \mid \mbox{$z \nin C$ 
and $\{z\} \cup C$ is a chain of $P$}\}.
\]

Clearly, $\tmop{lk}_P (C) = (< x_0) \ast (x_0, x_1) \ast \ldots \ast
(x_{n - 1}, x_n) \ast (> x_n)$ when $C \neq \emptyset$, while $\tmop{lk}_P (C) = P$ when $C = \emptyset$.

Then $P$ is said to be {\em Cohen-Macaulay} over the ring $\mathbb{A}$
  if for every chain $C$ in $P$, including $C = \emptyset$,
\[
\widetilde{H}_i (\tmop{lk}_P (C), \mathbb{A}) = 0
  \quad\mbox{for }  i \neq l (\tmop{lk}_P (C)) .
  \]

Unless stated otherwise, `Cohen-Macaulay' will mean `Cohen-Macaulay over $\mathbb{A}$'.

\begin{lemma}(Garst {\cite{Garst1979}})\label{lemma:garst}
The poset $P$ is Cohen-Macaulay over the ring
  $\mathbb{A}$ if and only if $\widetilde{H}_i (P,\mathbb{A}) = 0$ for $i \neq l (P)$, and $\tmop{lk}_P (x)$ is Cohen-Macaulay for all $x \in P$.
\end{lemma}

It is known (see \cite[Lemma 1.16]{Rylands1990}, \cite[Exercise 4.1.3]{Wachs2004}) that a
Cohen-Macaulay poset is ranked. The following proposition is also well-known (see \cite[Proposition 1.14]{Rylands1990} or \cite[Proposition 2.4.7]{Koonin2012} for a proof):

\begin{proposition}{\cite[Proposition 1.14]{Rylands1990}}\label{prop:CMjoin}
Let $P_1, \ldots, P_n$ be posets. Then $P_1
  \ast \cdots \ast P_n$ is Cohen-Macaulay over the field $\mathbb{F}$ if and
  only if $P_1, \ldots, P_n$ are Cohen-Macaulay over $\mathbb{F}$.
\end{proposition}

In particular,

\begin{corollary}
  \label{corollary:joinCM} Let $P^a = P \ast \{a\}$, where $\{a\}$ is the
  poset with one element. Then $P$ is Cohen-Macaulay over $\mathbb{F}$ if and
  only if $P^a$ is Cohen-Macaulay over $\mathbb{F}$. Similarly let $P_a
  =\{a\} \ast P.$ Then $P$ is Cohen-Macaulay over $\mathbb{F}$ if and only if
  $P_a$ is Cohen-Macaulay over $\mathbb{F}$.
\end{corollary}

\subsection{Unitary reflection groups and the Shephard-Todd classification}
\label{subsection:classification}

Let $V$ be a vector space over $\mathbb{C}$ of dimension $n$.

\begin{definition}
  An element $g \in GL(V)$ is a {\em reflection} if the
  order of $g$ is finite and if $\dim (\tmop{Fix} g) = n - 1$.  If $g$ is a
  reflection, the subspace $\tmop{Fix} g$ is a hyperplane, called the
  reflecting hyperplane of $g$. 
\end{definition}

\begin{definition}
  A {\em unitary reflection group} is a finite subgroup of $GL(V)$, which is generated
  by reflections.  
\end{definition}

Any such group is a group of unitary reflections with respect to some hermitian form, which explains the term unitary.

If $G$ is a unitary reflection group, the $G$-module
$V$ is called the {\em natural} (or {\em reflection})
representation of $G$. If $V$ is irreducible as a
$G$-module, we say that $G$ is an {\em irreducible
unitary reflection group}. The irreducible unitary reflection groups were
first classified by Shephard and Todd \cite{ShTo1954}.

There is an infinite family $G (r, p, n)$ of unitary reflection groups indexed
by triples of positive integers $(r, p, n)$ such that $p$ is a
divisor of $r$. The group $G (r, p, n)$ is defined to be the subgroup of $GL_{n} (\mathbb{C})$ consisting
of all $n \times n$ monomial matrices whose non-zero entries are complex $r$th
roots of unity such that the product of the non-zero entries is an
$\frac{r}{p}$th root of unity.

In addition, there are 34 {\em exceptional} irreducible unitary
reflection groups. This includes (complexifications of) the six exceptional
irreducible real reflection groups (of types $H_3, H_4, F_4, E_6, E_7
\tmop{and} E_8)$.

\subsection{Invariant theory for unitary reflection groups}\label{subsection:invariant}

Let $S$ denote the algebra of polynomial functions on $V$, which may be
identified with the symmetric algebra of the dual space, $S (V^{\ast})$.
Then $GL(V)$ acts on $S$, and if $G \subseteq GL(V)$ is a
finite group, we denote by $S^G$ the subalgebra of all $G$-invariant functions.
Both $S$ and $S^G$ have a natural grading by polynomial degree.

It is well-known that if $G$ is a unitary reflection group then $S^G$ is free,
and in fact this property characterises unitary reflection groups among finite
subgroups of $GL(V)$ (see \cite[Theorem 2.4]{Springer1974}).

Let $F$ denote the ideal of $S$ generated by the elements of $S^G$ with
nonzero constant term. The quotient $S / F$ is called the
{\em coinvariant algebra} of $G$ and is denoted $S_G$. The coinvariant
algebra $S_G$ inherits a grading from $S$. It is known (see \cite[Corollary 3.29 and Proposition 3.32]{LeTa2009} that $\dim S_G = \left| G \right|$ and that
the representation of $G$ on $S_G$ is the regular representation.

Let $t$ be an indeterminate, and $M$ any finite-dimensional $G$-module. Define the
{\em fake degree} of $M$ to be the polynomial
\[
f_M (t) \assign \sum_i \langle (S_G)_i, M \rangle_{}
t^i = \sum^r_{j = 1} t^{q_j (M)}
\]
where $\langle, \rangle$ denotes the usual intertwining number for complex representations of $G$. The integers $q_1 (M) \leqslant q_2 (M)
\leqslant \ldots \leqslant q_r (M)$ are called the $M$-exponents of $G$. Since $S_G$ affords the regular representation of $G$, $r = \dim M$.

In the case $M = V$ set $f_V (t) = \sum_j t^{m_j}$. The integers $m_i = q_i
(V)$ are called the {\em exponents} of $G.$ It is known (see \cite[Corollary 10.23]{LeTa2009}) that $m_i = d_i - 1$, where $d_{1}, \ldots, d_r$ are
the degrees of the basic invariants of $G$. 

When $M = V^{\ast}$, the corresponding exponents $m_i^{\ast} \assign q_i
(V^{\ast})$ are called the {\em coexponents} of $G$. Define
the {\em codegrees} $d_i^{\ast}$ of $G$ by $d_i^{\ast} \assign
m_i^{\ast} - 1$.

Let $\mathcal{H}$ be the space of $G$-harmonic polynomials -- that is,
those polynomials which are annihilated by all $G$-invariant
polynomial differential operators on $S$ with no constant term (see \cite[\S\S9.5-6]{LeTa2009} for more details). Clearly $\mathcal{H}$ is graded by degree. It is known (\cite[Corollary 9.37]{LeTa2009}) that $\mathcal{H}$ is a $G$-stable
complement of $F$ in $S$, so that $\mathcal{H}$ may be identified with $S_G$.
If $\mathcal{N}= N_{GL (V)} (G)$ then $\mathcal{H}$ is
$\mathcal{N}$-stable (as is $S^G$), and there is an isomorphism $S = S^G \otimes
\mathcal{H}$ of $\mathcal{N}$-modules (see \cite[Proposition 12.2]{LeTa2009}).

If $M$ is any finite-dimensional $G$-module, define a grading on the
$G$-module $S \otimes M$ by defining the degree $j$ component of $S
\otimes M$ to be $S_j \otimes M$. The module $(\mathcal{H} \otimes
M^{\ast})^G$ has a homogeneous linear basis $u_1, \ldots, u_r$, where $r =
\dim M$, and the degree of $u_i$ is $q_i (M)$. Such a basis is also an
$S^G$-basis of $(S \otimes M^{\ast})^G$ (see \cite[Proposition 10.3]{LeTa2009}). The elements $u_1, \ldots, u_n$ are called {\em covariants} of
$G$. 

\subsection{Parabolic subgroups}

\begin{definition}
  Suppose $G\subseteq GL(V)$ is a unitary reflection group on $V$, and that $U$ is
  a subset of $V$. Denote by $G_U$ the subgroup
\[
\{g \in G \mid \tmop{Fix} (g) \supseteq U\}.
\]
These groups are called {\em parabolic subgroups} of $G$. This generalises the notion of a parabolic subgroup for a real reflection group.
\end{definition}

The following is known as Steinberg's fixed point theorem (\cite[Theorem 1.5]{Steinberg1964}). For a more recent proof, see \cite{Lehrer2004}, and \cite[\S9.7]{LeTa2009}.

\begin{theorem}
  \label{theorem:steinberg}Let $G$ be a finite reflection group. If $v \in V$, the stabiliser $G_v =\{g \in G \mid gv =
  v\}$ is the reflection group generated by the reflections which fix v.
\end{theorem}

\begin{corollary}{\cite[Corollary 9.51]{LeTa2009}}
  \label{corollary:steinberg} Let $U$ be any
  subset of $V$. The parabolic subgroup $G_U$ of $G$ which fixes $U$ pointwise is
  the reflection group generated by the reflections in $G$ whose reflecting
  hyperplanes contain $U$. 
\end{corollary}

\subsection{Hyperplane arrangements and posets of
eigenspaces}\label{subsection:hyperplanes}

For background on hyperplane arrangements, see \cite{OrTe1992}.  For a discussion of hyperplane arrangements in the context of unitary reflection groups, see \cite{OrSo1982}.  In the following section, let
$V$ be a complex vector space of dimension $n$.

A {\em hyperplane} of $V$ is a subspace $H$
of codimension 1 in $V$ - that is, a linear hyperplane.  A {\em hyperplane arrangement}
$\mathcal{A}$ in $V$ is a finite set of hyperplanes in $V$. A {\em subarrangement} $\mathcal{A}'$ of $\mathcal{A}$ is a subset of
$\mathcal{A}$.  A hyperplane arrangement is
{\em essential} if $\cap_{H \in \mathcal{A}} H =\{0\}$. 

\begin{definition}
 Let $G$ be a unitary reflection group, and let
$\mathcal{A}=\mathcal{A}(G)$ be the set of reflecting hyperplanes of
$G$ -- i.e. the set of hyperplanes of $V$ which are the fixed
point subspaces of elements in $G$. Define $\mathcal{L}(\mathcal{A}(G))$ to be the poset of
intersections of the hyperplanes in $\mathcal{A}$. By convention, $V \in
\mathcal{L}(\mathcal{A}(G))$ -- it is the empty intersection of elements of $\mathcal{L}(\mathcal{A}(G))$. 
\end{definition}

It is well known that $\mathcal{L}(\mathcal{A}(G))$ is a geometric lattice (see \cite[Lemma 2.3]{OrTe1992}), and hence (see \cite[Theorem 4.1]{Folkman1966}) that it has
vanishing reduced homology except in top dimension. 

\begin{theorem}{\cite [Theorem 6.27]{OrTe1992}}\label{proposition:fixedpoint} Let $G \subseteq GL(V)$ be a unitary reflection group. Then:
\begin{itemize}
\item[(i)] If $g \in G$, then $\tmop{Fix} (g) \in \mathcal{L}$.
\item[(ii)] If $X \in \mathcal{L}(\mathcal{A}(G))$, then there exists $g \in G$ with $\tmop{Fix}
  (g) = X$.
  \end{itemize}
\end{theorem}

Thus the set of subspaces $\{\tmop{Fix} g \mid g \in G\}$ coincides with $\mathcal{L}(\mathcal{A}(G))$.

\subsection{Reflection cosets}\label{subsection:coset}
  
\begin{definition}
  Suppose $G$ is a unitary reflection group in $V$. Let
  $\gamma$ be a linear transformation of finite order in $V$, such that $\gamma
  G = G \gamma$. The coset $\gamma G$ is called a {\em reflection coset}.
\end{definition}

Write $\widetilde{G}=\langle G, \gamma \rangle$ for the group generated by $G$ and $\gamma$. Note that $\widetilde{G}$ is finite.

Let $M$ be any finite-dimensional $G$-module. Recall from
\S\ref{subsection:invariant} that the module $(\mathcal{H} \otimes
M^{\ast})^G$ has a homogeneous linear basis $u_1, \ldots, u_r$, where $r =
\dim M$, and the degree of $u_i$ is $q_i (M)$. In fact (see \cite[Proposition 12.2(ii)]{LeTa2009}) such a basis may be chosen such that for each $i$, $\gamma u_i = \varepsilon_i (M) u_i$ for some complex root of unity $\varepsilon_i (M) \in \mathbb{C}^{\times}$. Furthermore, the
(multi)set of pairs $(\deg u_i, \varepsilon_i (M))$ depends only on the coset
$\gamma G$ and the module $M$, and not on the particular choice of
$\gamma$ \cite[Proposition 12.2(iii)]{LeTa2009}. The multiset $\{\varepsilon_i
(M)\}$ is called the {\em set of $M$-factors of $\gamma G$}. In the
special case $M = V$, the $V$-factors $\varepsilon_i (V)$ will be denoted
$\varepsilon_i$ and called {\em factors}. Analogously, in the special
case $M = V^{\ast}$, the $V^{\ast}$-factors $\varepsilon_i (V^{\ast})$ will be
denoted $\varepsilon_i^{\ast}$ and called {\em cofactors.}

Irreducible reflection cosets were discussed by Cohen in \cite{Cohen1976} and classified in \cite{BrMaMi1999}, where they are referred to as {\em reflection data}. Also see \cite[Table D.5,
p.278]{LeTa2009}.

Let $\gamma G$ be a reflection coset in $V$. For $\zeta \in
\mathbb{C}^{\times}$ a root of unity, define
\[
A (\zeta) =\{i \mid
\varepsilon_i \zeta^{d_i} = 1\},
\]
where $d_i$ is the degree of the basic invariant corresponding to
$\varepsilon_i .$ Let $a (\zeta) = \left| A (\zeta) \right|$.  

The following are two key theorems of Springer-Lehrer theory:

\begin{theorem}{\cite[Theorem 12.19]{LeTa2009}} \label{theorem:eigenspacescoset} Let $\gamma
  G$ be a reflection coset in $V$, and $\zeta \in
  \mathbb{C}^{\times}$ be a complex root of unity. Suppose $E$ is maximal
  among $\zeta$-eigenspaces corresponding to elements of $\gamma G$. Then:
\begin{itemize}
\item[(i)] $\dim E = a (\zeta)$,
\item[(ii)] any two maximal eigenspaces $V (x, \gamma)$ are conjugate under the
  action of $G$. 
  \end{itemize}
\end{theorem}

\begin{theorem}{\cite[Theorem 12.20]{LeTa2009}}\label{theorem:subquotientcoset} Let $\gamma G$
  be a reflection coset in $V$, and for $\zeta \in \mathbb{C}^{\times}$, let $E
  = V (x, \zeta)$ be a $\zeta$-eigenspace, maximal among the
  $\zeta$-eigenspaces of elements of $\gamma G$. Let $N(E) =\{g \in G \mid gE =
  E\}$, and $C(E) =\{g \in G \mid gv = v \mbox{ for all } v \in E\}$. Then $N(E) / C(E)$ acts as a unitary reflection group on $E$, whose hyperplanes are the
  intersections with $E$ of the hyperplanes of $G$.
\end{theorem}

There is also a theory of regular elements for reflection cosets. As with reflection groups, the element $x \in \gamma G$ is $\zeta$-regular if $V (x,
\zeta)$ contains a $G$-regular vector. If $x \in \gamma G$ is
$\zeta$-regular for some $\zeta$, then $x$ is said to be regular. The following results will be important later:

\begin{proposition}{\cite[Proposition 12.21]{LeTa2009}} \label{proposition:coseteigenvalues} Let
  $\gamma G$ be a reflection coset and let $M$ be a $\widetilde{G}$-module of
  dimension $r$, where $\widetilde{G} = \langle G, \gamma \rangle$. If $x$ is a regular element
  of $\gamma G$, then the eigenvalues of $x$ on $M$ are $\{\varepsilon_i
  (M^{\ast}) \zeta^{q_i (M^{\ast})} \mid i = 1, \ldots, r\}.$
\end{proposition}

In particular,
  
  \begin{corollary}
    \label{corollary:coseteigenvalues} Let $M = V$, with the notation as in
    Proposition \ref{proposition:coseteigenvalues}.   If $x$ is a regular element
  of $\gamma G$, then the eigenvalues of
    $x$ acting on $V$ are $\{\varepsilon_i^{\ast} \zeta^{m_i^{\ast}}
    \mid i = 1, \ldots r\}$, where the $m_i^{\ast}$ are the coexponents of $G$,
    and the $\varepsilon_i^{\ast}$ are the cofactors defined above. 
  \end{corollary}

For a reflection group it is trivial that the element $1 \in G$ is regular. For a reflection coset the existence of a regular element is not obvious, and
indeed if $G$ acts reducibly on $V$, $\gamma$ could be chosen so that there
are no $\gamma G$-regular vectors in $V$.  See \cite[Remark 12.24]{LeTa2009}.  However, the following was proved in \cite[Corollary 7.3]{BoLeMi2006}, and independently in \cite[Theorem 3.4]{Malle2006}:

\begin{theorem}\label{theorem:cosetregular}
  If $\gamma G$ is a reflection coset and $G$ acts
  irreducibly on $V$, then there is a $\zeta$-regular element in $\gamma G$ for
  some $\zeta \in \mathbb{C}^{\times}$.
\end{theorem}

\subsection{Posets of eigenspaces}\label{subsection:posetsofeigenspaces}

The following definition contains the central object of study in this paper:

\begin{definition}
  \label{definition:S_v(zeta)} Let $\gamma G$ be a reflection coset in $V
  =\mathbb{C}^n$, and $\zeta \in \mathbb{C}^{\times}$ be a complex root of
  unity. Define $\mathcal{S}_{\zeta}^V (\gamma G)$ to be the set $\{V (x,
  \zeta) \mid x \in \gamma G\}$, partially ordered by the reverse of
  inclusion. 
\end{definition}

It will be seen (Corollary \ref{corollary:S_v(zeta)maximal}) that the poset
$\mathcal{S}_{\zeta}^V (\gamma G)$ always has a unique maximal element $\hat{1}$, and
that it may or may not have a unique minimal element $\hat{0}$ as
well (for example, the full space $V$). Thus the associated simplicial complex is a cone, which is contractible. Hence its homology is uninteresting, and we have reason to study the following poset:

\begin{definition}
  \label{definition:Stilde_v(zeta)} Define $\widetilde{\mathcal{S}}_{\zeta}^V (\gamma
  G)$ to be the subposet of $\mathcal{S}_{\zeta}^V (\gamma G)$ obtained by
  removing the unique maximal element, as well as the unique minimal element if it
  exists.
\end{definition}

\section{General structure theorems for $\mathcal{S}_V^{\zeta} (\gamma G)$}

\label{section:general}This section examines the structure of the posets $\mathcal{S}_{\zeta}^V
(\gamma G)$ for arbitrary unitary reflection groups $G$ and $\gamma
\in N_{GL(V)} (G)$. With the idea of showing that the posets
$\mathcal{S}_{\zeta}^V (\gamma G)$ are Cohen-Macaulay in mind, we consider
closed intervals in $\mathcal{S}_{\zeta}^V (\gamma G)$.  The key is the
following theorem:

\begin{theorem}\label{theorem:poset}
Let $\gamma G$ be a
  reflection coset in $V$. Let $E = V (\gamma g, \zeta) \in
  \mathcal{S}_{\zeta}^V (\gamma G)$. Then the following three posets are
  identical (as sets of subspaces of $V)$:
\begin{itemize}
\item[(i)] $\{V (\gamma h, \zeta) \in \mathcal{S}_{\zeta}^V (\gamma G) \mid h \in
  G, V (\gamma h, \zeta) \subseteq E\};$
\item[(ii)] $\{V (\gamma h, \zeta) \cap E \mid h \in G\};$
\item[(iii)] $\{E \cap \tmop{Fix}_V (x) \mid x \in G\}$.
\end{itemize} 
 In addition, if $E$ is maximal among the $\zeta$-eigenspaces of elements of
  $\gamma G$, then there is a fourth poset
\begin{itemize}
\item[(iv)] $\mathcal{S}_1^E (N (E) / C (E))$    
\end{itemize}
 equal to each of the above.
 \end{theorem}

Before embarking on the proof of this theorem, we point out some
  important corollaries:
  
\begin{corollary}
    \label{corollary:joinsemilattice}
    With notation as above, suppose $E_1, E_2 \in \mathcal{S}_{\zeta}^V (\gamma G)$.  Then $E_1 \cap E_2 \in \mathcal{S}_{\zeta}^V (\gamma G)$.  In particular, $\mathcal{S}_{\zeta}^V (\gamma G)$ is a join semilattice. 
  \end{corollary}
  
  \begin{proof}[{\em Proof of Corollary \ref{corollary:joinsemilattice}}]
    Assuming that the posets (i) and (ii) are equal, it follows that the
    intersection of any two $\zeta$-eigenspaces is again a $\zeta$-eigenspace.
    Hence any two elements of $\mathcal{S}_{\zeta}^V (\gamma G)$ have a
    least upper bound, and so $\mathcal{S}_{\zeta}^V (\gamma G)$ is a
    join semilattice.
  \end{proof}
  
  \begin{corollary}
    \label{corollary:S_v(zeta)maximal}With notation as above, the poset
    $\mathcal{S}_{\zeta}^V (\gamma G)$ always has a unique maximal element $\hat{1}$.
      \end{corollary}

    \begin{proof}[{\em Proof of Corollary \ref{corollary:S_v(zeta)maximal}}]
      By repeated application of Corollary \ref{corollary:joinsemilattice}, $\bigcap_{E \in
      \mathcal{S}_{\zeta}^V (\gamma G)} E$ is an eigenspace, and by
      definition it is contained within every element $E \in
      \mathcal{S}_{\zeta}^V (\gamma G)$.
    \end{proof}
  
  \begin{remark}
    As a result of Corollary \ref{corollary:S_v(zeta)maximal}, the first poset mentioned in Theorem \ref{theorem:poset} can be
    written $[E, \hat{1}] .$
  \end{remark}
  
  \begin{corollary}
    Suppose $E = V (\gamma g, \zeta) \in \mathcal{S}_{\zeta}^V (\gamma G)$.
    \label{corollary:S_v(zeta)upper} The interval $[E, \hat{1}]$ is a
    geometric lattice.
  \end{corollary}
  
  \begin{proof}[{\em Proof of Corollary \ref{corollary:S_v(zeta)upper}}]
    The eigenspace $E = V (\gamma g, \zeta)$ is contained in some maximal
    eigenspace $E'$. Hence by Theorem \ref{theorem:poset}(iv), $[E, \hat{1}]$
    is a principal upper order ideal of $\mathcal{S}_1^{E'} (N (E') / C (E'))
    .$ The latter is known to be a geometric lattice by Proposition \ref{proposition:fixedpoint} and \cite[Lemma 2.3]{OrTe1992}. Since closed intervals of geometric
    lattices are also geometric lattices (see \cite[Lecture 3]{Stanley2007}), the
    result follows.
  \end{proof}
  
We now prove Theorem \ref{theorem:poset}.
  
\begin{proof}[{\em Proof of Theorem \ref{theorem:poset}}]

The equality of (iii) and (iv) follows easily from Theorem  \ref{theorem:subquotientcoset}.

Hence we concentrate on the equality of (i), (ii) and (iii).

 (ii) = (iii). First we show that, for any element $h \in G$, $V (\gamma
    h, \zeta) \cap E = E \cap \tmop{Fix}_V (\gamma g h^{- 1} \gamma^{- 1})$.
    Note that $\gamma g h^{- 1} \gamma^{- 1} \in G$ since $\gamma$ normalises
    $G.$
    
	If $v \in V (\gamma h, \zeta) \cap E$ then $\gamma gv = \zeta v$ and
    $\gamma hv = \zeta v$. So $v = (\gamma g) (\gamma h)^{- 1} v =
    (\gamma g h^{- 1} \gamma^{- 1}) v$. Hence $v \in \tmop{Fix}_V (\gamma g
    h^{- 1} \gamma^{- 1})$. Conversely, if $v \in E \cap \tmop{Fix}_V (\gamma
    g h^{- 1} \gamma^{- 1})$ then $\gamma g v = \zeta v$, and $\gamma g h^{-
    1} \gamma^{- 1} v = v.$ Hence $(\gamma h)^{- 1} v = h^{- 1} \gamma^{- 1}
    v = (\gamma g)^{- 1} (\gamma g h^{- 1} \gamma^{- 1}) v = (\gamma g)^{- 1}
    v = \zeta^{- 1} v.$ Thus $(\gamma h) v = \zeta v, \tmop{and} \tmop{so} v
    \in V (\gamma h, \zeta)$. This proves that $V (\gamma h, \zeta) \cap E =
    E \cap \tmop{Fix}_V (\gamma g h^{- 1} \gamma^{- 1})$, and therefore that
    (ii)$\;\subseteq\;$(iii). To prove the reverse inclusion, if $x$ is
    any element of $G$, let $h = \gamma^{- 1} x^{- 1} \gamma g$. \
    Substituting into the formula just proved, we have $V (\gamma h, \zeta)
    \cap E = E \cap \tmop{Fix}_V (x)$. Thus the reverse inclusion is
    true, and (ii) = (iii).

    By definition, (i)$\;\subseteq\;$(ii). Hence the theorem will be
    proved if we can show that (iii)$\;\subseteq \;$(i). In other words, we must
    prove
    
    \begin{proposition}\label{prop:int1}
 Let $V$ be a complex vector space, and
      $\gamma G$ a reflection coset in $V$. Suppose that $E = V (\gamma, \zeta) \in \mathcal{S}_{\zeta}^V (\gamma G)$, and $X \in \mathcal{L}(G)$.  Then $E
      \cap X = V (\gamma h, \zeta)$ for some $h \in G$.
         \end{proposition}

 Note that we can
      assume without loss of generality that $E = V (\gamma, \zeta)$ rather than $V (\gamma g, \zeta)$ by replacing $\gamma$
      with $\gamma g$, since if $\gamma \in N_{GL(V)} (G)$ then $\gamma g
\in N_{GL(V)} (G)$ also.
    
Consider the following three propositions:
    
    \begin{proposition}\label{prop:int2}
Let $V$ be a complex vector space, and
      $\gamma G$ any reflection coset in $V.$  Suppose that $E
      = V (\gamma, \zeta) \in \mathcal{S}_{\zeta}^V (\gamma G)$ and $X \in \mathcal{L}(G)$ such that $\gamma X = X$. Then $E \cap X = V (\gamma h, \zeta)$ for some $h \in G$. 
    \end{proposition}
    Note that we are not
      assuming that $E$ is maximal.
    
    \begin{proposition}\label{prop:int3} 
Let $V$ be a complex vector space, and
      $\gamma G$ any reflection coset in $V$, where $G$ acts
      essentially on $V$ (i.e. $\cap_{H \in \mathcal{A}_G} H =\{0\}$). Then
      there exists $h \in G$ such that $V (\gamma h, \zeta) =\{0\}.$
    \end{proposition}
    
    \begin{proposition}\label{prop:int4}
Let $V$ be a complex vector space, and
      $ \gamma G$ any reflection coset in $V$, where $G$ acts
      essentially on $V$ (i.e. $\cap_{H \in \mathcal{A}_G} H =\{0\})$. Then
      there exists $h \in G$ such that $V (\gamma h, 1) =\{0\}.$
    \end{proposition}
    
    The key to proving Theorem \ref{theorem:poset} is the following reduction
    theorem:
    
    \begin{theorem}      \label{theorem:implications1}
We have the implications
    (\ref{prop:int4}) $\Rightarrow$ (\ref{prop:int3}) $\Rightarrow
      $ (\ref{prop:int2}) $\Rightarrow$ (\ref{prop:int1}).
    \end{theorem}
    
    \begin{proof}[{\em Proof of Theorem \ref{theorem:implications1}}]
      (\ref{prop:int2}) $\Rightarrow$ (\ref{prop:int1}). Let $E$ and $X$ be
      as in Proposition \ref{prop:int1}. Define $X' \assign \bigcap_{i =
      0}^{\infty} \gamma^i X.$ Note that this intersection is actually finite
      as $\gamma$ has finite order. For each $i$, $\gamma^i X \in \mathcal{L}(G)$, since $\gamma$ normalises
      $G$. Thus $X' \in \mathcal{L}(G)$, and so by Proposition
      \ref{proposition:fixedpoint}(ii), $X' = \tmop{Fix}_V (y)$ for some $y \in
      G.$

      Clearly $\gamma X' = X' .$ Applying Proposition \ref{prop:int2} to
      $X'$, $E \cap X' = V (\gamma h, \zeta)$ for some $h \in G.$

      Now $E \cap X' = E \cap (\cap_{i = 0}^{\infty} \gamma^i X) = \cap_{i =
      0}^{\infty} (E \cap \gamma^i X)$. Also, 
      \begin{align*}
      E \cap \gamma^i X &= \gamma^{- i} (E \cap \gamma^i X)\quad\mbox{ (since
      $\gamma$ acts by scalar multiplication on $E$)}\\
 &= \gamma^{- i} E \cap X\\
&= E \cap X.
\end{align*}
Hence $V (\gamma h, \zeta) = E \cap X' = \bigcap_{i = 1}^{\infty} (E \cap
      \gamma^i X) = \bigcap_{i = 1}^{\infty} (E \cap X) = E \cap X$, as required.

      (\ref{prop:int3}) $\Rightarrow$ (\ref{prop:int2}). Let $E$ and $X$ be
      as in Proposition \ref{prop:int2}. Now $G_X$ acts
      as a reflection group on $V$ by Steinberg's fixed point theorem
      (Corollary \ref{corollary:steinberg}), and it acts essentially on $V / X$.       Since $\gamma X = X$ we have $\gamma \in N_{GL(V)} (G_X)$.     Hence $\gamma G_X$ is a reflection coset in $V / X$.

      Applying Proposition \ref{prop:int3} to $\gamma G_X,$ there is some $h
      \in G_X$ such that $(V / X) (\gamma h, \zeta) =\{0\}.$ In other words,
      $V (\gamma h, \zeta) \subseteq X$. Also if $v \in V (\gamma h, \zeta)$
      and $h \in G_X$, then
\[
\gamma v = \gamma (h v) = (\gamma h) v = \zeta v.
\]

      Hence $V (\gamma h, \zeta) \subseteq E.$ Thus $V (\gamma h, \zeta)
      \subseteq E \cap X$. Conversely, $E \cap X \subseteq V (\gamma h,
      \zeta)$, so $E \cap X = V (\gamma h, \zeta)$.

      (\ref{prop:int4}) $\Rightarrow$ (\ref{prop:int3}). If $\gamma G$ is a
      reflection coset in $V$, then so is $\zeta^{- 1} \gamma G$. Applying
      Proposition \ref{prop:int4} to $\zeta^{- 1} \gamma G$, there exists $h
      \in G$ such that $V (\zeta^{- 1} \gamma h, 1) =\{0\}$. But $V (\gamma
      h, \zeta) = V (\zeta^{- 1} \gamma h, 1)$, as required.
    \end{proof}
    
    We now focus on the proof of Proposition \ref{prop:int4}. The first task
    is to reduce Proposition \ref{prop:int4} to the case where $G$ is an
    irreducible reflection group. Consider the following two Propositions:
    
    \begin{proposition}\label{prop:int5}
      Proposition \ref{prop:int4} holds when $V$ is
      irreducible as a $\widetilde{G}$-module, where $\widetilde{G} = \langle G, \gamma \rangle$.
    \end{proposition}
    
    \begin{proposition}\label{prop:int6}
  Proposition \ref{prop:int4} holds when $V$ is
      irreducible as a $G$-module.
    \end{proposition}
    
    Clearly (\ref{prop:int4}) $\Rightarrow$ (\ref{prop:int5}) $\Rightarrow$ (\ref{prop:int6}). The reverse implications are also true:
    
    \begin{theorem}\label{theorem:implications2}
  We have the implications
      (\ref{prop:int6}) $\Rightarrow$ (\ref{prop:int5}) $\Rightarrow$ (\ref{prop:int4}).
    \end{theorem}
    
    \begin{proof}[{\em Proof of Theorem \ref{theorem:implications2}}]
      (\ref{prop:int5}) $\Rightarrow$ (\ref{prop:int4}). Suppose $V =
      \bigoplus_{i = 1}^k V_i$ is a decomposition of $V$ into irreducible
      $\widetilde{G}$-modules. Then $G = G_1 \times G_2 \times \ldots \times G_k$
      and $\gamma = \gamma_1 \oplus \gamma_2 \oplus \ldots \ldots \oplus
      \gamma_k,$ where $G_i$ acts on $V_i$ as a reflection group and trivially
      on $V_j$ for $j \neq i$, and $\gamma_i$ normalises $G_i$.

      Note that $\mathcal{A}(G) = \cup_{i = 1}^k \mathcal{A}(G_i)$.

      Thus $\bigcap_{H \in \cup \mathcal{A}(G_j), j \neq i} H \supseteq V_i$
      (and by hypothesis equals $V_i) .$

      It follows that $\bigcap_{H \in \mathcal{A}(G_i)} H = \bigoplus_{j \neq
      i} V_j$, and hence that $\bigcap_{H \in \mathcal{A}(G_i)} (H \cap V_i)
      =\{0\}.$

      Hence $(G_i, \gamma_i, V_i)$ satisfies the conditions of Proposition
      \ref{prop:int5}. So by assumption there exists elements $x_i \in G_i$
      such that $V_i (\gamma_i x_i, 1) =\{0\}.$

      Hence $V (\gamma x, 1) =\{0\},$ where $x = x_1 \oplus x_2 \oplus \ldots
      \oplus x_k,$ and Proposition \ref{prop:int4} follows.

      (\ref{prop:int6}) $\Rightarrow$ (\ref{prop:int5}). Assume that $V$ is
      irreducible as a $\widetilde{G}$-module. Then $G = G_1 \times G_2 \times
      \ldots \times G_k$, $V = V_1 \oplus V_2 \oplus \ldots \oplus V_k$, and all
      $(G_i, V_i)$ are isomorphic, irreducible, and permuted cyclically by
      $\gamma$. Thus $\gamma^k$ fixes all the $V_i$ and in particular
      normalises $G_i$ for each $i$. (See \cite[Proposition 6.9]{BlLe2001} for further
      details.)

      This means that $V_i = \gamma^{i - 1} V_1$,  $i = 1, 2, \ldots, k$ and
      $\gamma^k V_1 = V_1$. Thus we may write an arbitrary element $v \in V$
      uniquely as $v = v_1 + \gamma v_2 + \ldots + \gamma^{k - 1} v_k$ where
      $v_i \in V_1$ for all $i$.

      Then 
      \begin{align*}
      \gamma v &= \gamma ( v_1 \oplus \gamma v_2 \oplus \ldots \oplus
      \gamma^{k - 1} v_k)\\
      &= \gamma^k v_k \oplus \gamma v_1 \oplus \ldots
      \oplus \gamma^{k - 1} v_{k - 1}.
      \end{align*}
 Thus $v$ is fixed by $\gamma$ if and only if $v_1 = v_2 = \ldots = v_k
      \tmop{and} \gamma^k v_1 = v_1 .$

      Now by Proposition \ref{prop:int6} applied to $(V_1, \gamma^k, G_1)$,
      there exists $y \in G_1$ such that $\tmop{Fix}_{V_1} (\gamma^k y) = V_1
      (\gamma^k y, 1) =\{0\}$.

      We may write $x \in G$ uniquely as $x = (x_1, \text{}^{\gamma} x_2,
      \ldots, \text{}^{\gamma^{k - 1}} x_k) \text{}^{}$, where $x_i \in G_1$
      for all $i$, and the left superscript denotes conjugation.

      Then 
      \begin{align*}
      (\gamma x) v &= \gamma (x v)\\
      &=\gamma (x_1 v_1 \oplus \gamma x_2 v_2
      \oplus \ldots \oplus \gamma^{k - 1} x_k v_k)\\
      &=\gamma^k x_k v_k \oplus \gamma x_1 v_1
      \oplus \ldots \oplus \gamma^{k - 1} x_{k - 1} v_{k - 1}.
      \end{align*}
      
This equals $v = v_1 \oplus \gamma v_2 \oplus \ldots \oplus \gamma^{k -
      1} v_k$ if and only if
\[
\gamma^k x_k v_k = v_1, x_1 v_1 =
      v_2, x_2 v_2 = v_3, \ldots, x_{k - 1} v_{k - 1} = v_k.
      \]
That is,
\[
v_k = x_{k - 1} v_{k - 1} = x_{k - 1} x_{k - 2} v_{k - 2} =
      \ldots = x_{k - 1} x_{k - 2} \ldots x_2 v_2 = x_{k - 1} x_{k - 2} \ldots
      x_1 v_1.
      \]

      Therefore $\gamma^k x_k \ldots x_1 v_1 = v_1$.

      To summarise, we have proved that if $v \in \tmop{Fix}_V (\gamma x) = V
      (\gamma x, 1)$, then $\gamma^k x_k x_{k - 1} \ldots x_1 v_1 = v_1
      .$

      Now take $x = (y, 1, 1, \ldots, 1)$. Then if $(\gamma x) v = v$, we
      must have $(\gamma^k y) v_1 = v_1$, whence by assumption $v_1 = 0.$ \
      But then $v_i = 0$ for all $i$, and so $v = 0$. Hence $V
      (\gamma x, 1) =\{0\}$, as required.
    \end{proof}
          
      Therefore, in view of Theorems \ref{theorem:implications1} and
      \ref{theorem:implications2}, Theorem \ref{theorem:poset} will follow
      from Proposition \ref {prop:int6}.
 
Before proving Proposition \ref {prop:int6}, we record a lemma which will be useful later:

\begin{lemma}\label{lemma:Ggammairreducible}
With notation as above, $V (\gamma x, \zeta)$ and 
$V (\gamma^k x_k x_{k-1} \ldots x_1, \zeta^k)$ are isomorphic as vector spaces.
\end{lemma}

\begin{proof}
As in the proof of Theorem \ref{theorem:implications2}, 

      \begin{align*}
      (\gamma x) v &= \gamma (x v)\\
      &=\gamma (x_1 v_1 \oplus \gamma x_2 v_2
      \oplus \ldots \oplus \gamma^{k - 1} x_k v_k)\\
      &=\gamma^k x_k v_k \oplus \gamma x_1 v_1
      \oplus \ldots \oplus \gamma^{k - 1} x_{k - 1} v_{k - 1}.
      \end{align*}
      
This equals $\zeta v = \zeta  v_1 \oplus \zeta \gamma v_2 \oplus \ldots \oplus \zeta \gamma^{k -
      1} v_k$ if and only if
\[
\gamma^k x_k v_k = \zeta v_1, x_1 v_1 = \zeta
      v_2, x_2 v_2 = \zeta v_3, \ldots, x_{k - 1} v_{k - 1} = \zeta v_k.
      \]
That is,
\[
v_k = \zeta^{-1} x_{k - 1} v_{k - 1} = \zeta^{-2} x_{k - 1} x_{k - 2} v_{k - 2} =
      \ldots = \zeta^{-(k-2)} x_{k - 1} x_{k - 2} \ldots x_2 v_2 = \zeta^{-(k-1)} x_{k - 1} x_{k - 2} \ldots
      x_1 v_1.
      \]

      Therefore $\gamma^k x_k \ldots x_1 v_1 = \zeta^k v_1$. That is, $v_1 \in V (\gamma^k x_k x_{k-1} \ldots x_1, \zeta^k)$. On the other hand, the above equations make it clear that if $v = v_1 + \gamma v_2 + \ldots + \gamma^{k - 1} v_k \in V (\gamma x, \zeta)$ then $v_2, \ldots, v_k$ are uniquely determined by $v_1$.  Thus the map which sends $v$ to $v_1$ is a vector space isomorphism.
\end{proof}

 We now prove Proposition \ref{prop:int6}.
        
        \begin{proof}[{\em Proof of Proposition \ref{prop:int6}}]
          The proof starts with a further reduction:
          
          \begin{lemma}\label{lemma:innerautomorphism}
            If $\gamma$ induces an
            inner automorphism of $G$, then Proposition
            \ref{prop:int6} holds.
            \end{lemma}
            
            \begin{proof}[{\em Proof of Lemma \ref{lemma:innerautomorphism}}]
              By hypothesis, there exists $h \in G$ such that $\gamma h$ is
              central in $\langle G, \gamma \rangle$ and hence is scalar. If $\gamma h \neq 1$ then $\tmop{Fix}_V
              (\gamma h) =\{0\}$, while if $\gamma h = 1$, then $\gamma
              \in G$ and $\gamma G = G$. In the latter case, since $G
              \neq \{1\}$, there exists at least one reflecting hyperplane and
              so $\{0\}= \cap_{\alpha} H_a = \tmop{Fix}_V (g)$ for some $g \in G$ (by Proposition
              \ref{proposition:fixedpoint}).
            \end{proof}

            In \cite[Table D.5, p.278]{LeTa2009}, there is a list of reflection
            cosets $(G, \gamma_0)$. Any reflection coset $(G, \gamma)$
            satisfying the conditions of Proposition \ref{prop:int6}, in which
            $\gamma$ does not induce an inner automorphism, has the form
            $(G, \gamma) = (G, \alpha \gamma_0)$, where $(G,
            \gamma_0)$ appears in this table, and $\alpha \in
            \mathbb{C}^{\times}$. (See also \cite[Proposition 3.13]{BrMaMi1999}).             We check Proposition \ref{prop:int6} for all such cosets case by
            case. That is, it is necessary to show that the reflection coset
            $(\alpha \gamma_0) G$ contains an element $x$ such that
            $\tmop{Fix}_V (x) =\{0\}$, where $\alpha \in
            \mathbb{C}^{\times}$. \

            The proof relies on the following lemma:
            
            \begin{lemma}\label{lemma:eigenvalue}
              Suppose $\alpha \gamma_0 G$ is a
              reflection coset in V which satisfies the conditions of
              Proposition \ref{prop:int6}. If there exists $x \in \alpha
              \gamma_0 G$ which acts on $V$ in such a way that none of its
              eigenvalues is equal to 1, then $\tmop{Fix}_V (x) =\{0\}$.
            \end{lemma}
            
            \begin{proof}[{\em Proof of Lemma \ref{lemma:eigenvalue}}]
              This is immediate from the definition of an eigenvalue.
            \end{proof}
            
            We now proceed to check that Lemma \ref{lemma:eigenvalue} holds
            for all reflection cosets satisfying the conditions of
            Proposition \ref{prop:int6}. Let $\xi_m$ denote a fixed
            primitive $m$-th root of unity.
                        
           {\bf Case 1}: $G = G (m, p, n), \gamma = \tmop{diag}
            (\xi_{\frac{em}{p}}, 1, \ldots, 1)$, where $e \mid p$ and $\xi_{\frac{em}{p}}$ is a fixed primitive
            $\frac{em}{p}$-th root of unity. Note that \cite[Table D.5, p.278]{LeTa2009} treats the cases $m \neq p$ and $m = p$ separately.
            For the purposes of this argument they can be treated together.

            In this case $\alpha \gamma = \tmop{diag} (\alpha
            \xi_{\frac{em}{p}}, \alpha, \ldots, \alpha)$. The eigenvalues of
            $\alpha \gamma$ acting on $V$ are $\alpha \xi_{\frac{em}{p}}
            \tmop{and} \alpha .$ Hence if $\alpha \neq 1$,
            $\xi_{\frac{em}{p}}^{- 1}$ by Lemma \ref{lemma:eigenvalue} there is
            nothing more to prove.

            If $\alpha = 1$ we may assume that $e \neq 1$ (and hence $p, m
            \neq 1$, as $e \mid p)$, since if $e = 1,$ $\alpha \gamma$
            induces an inner automorphism. Thus $\xi_{\frac{em}{p}} \neq 1.$

            If $n$ is even define $y \in G$ by $y_{1, 2} = y_{2, 1} = 1$, $y_{i,
            i} = (\xi_m)^{(- 1)^{i + 1}}$ if $i \geqslant 3$, $y_{i, j} = 0$ otherwise. Then $y \in G$ and $x \assign \alpha
            \gamma y$ is defined by $x_{1, 2} = \xi_{\frac{em}{p}}, x_{2, 1} =
            1, x_{i, i} = (\xi_m)^{(- 1)^{i + 1}}$ if $i \geqslant 3$,
            $x_{i, j} = 0$ otherwise. It is easy to check that
            $\tmop{Fix}_V (x) =\{0\}.$

            If $n$ is odd define $y \in G$ by $y \assign \tmop{diag} (1,
            \xi_m, \xi_m^{- 1}, \ldots, \xi_m, \xi_m^{- 1})$. Then $x
            \assign \alpha \gamma y$ is equal to $\tmop{diag}
            (\xi_{\frac{em}{p}}, \xi_m, \xi_m^{- 1}, \ldots, \xi_m, \xi_m^{-
            1})$ and $\tmop{Fix}_V (x) =\{0\}.$

            If on the other hand $\alpha = \xi_{\frac{em}{p}}^{- 1}$ then
            $\alpha \gamma = \tmop{diag} (1, \xi_{\frac{em}{p}}^{- 1}, \ldots,
            \xi_{\frac{em}{p}}^{- 1})$. Again we may assume that $e \neq 1$
            as otherwise $\alpha \gamma$ induces an inner automorphism.
            Define $y \in G$ by $y_{1, 2} = y_{2, 1} = 1, y_{i, i} = 1$ if $i \geqslant 3$, $y_{i, j} = 0$ otherwise. Then
            $x \in \alpha \gamma G$ is defined by $x_{1, 2} =
            \xi_{\frac{em}{p}}^{- 1}, x_{2, 1} = 1, x_{i, i} =
            \xi_{\frac{em}{p}}^{- 1}$ if $i \geqslant 3$, $x_{i, j} = 0$ otherwise. Once again, it is easy to check that
            $\tmop{Fix}_V (x) =\{0\}.$ This completes the proof of Lemma
            \ref{lemma:eigenvalue} for Case 1.

            The proof of Proposition \ref{prop:int6} in the remaining cases relies on Corollary \ref{corollary:coseteigenvalues}, Theorem
            \ref{theorem:cosetregular}, and the following lemma:
            
   \begin{lemma}\label{lemma:cosetalpha}
   If $\gamma G$ is a
                reflection coset, $\alpha$, $\zeta \in \mathbb{C}^{\times}$,
                and $x \in \gamma G$ then:
\begin{itemize}
\item[(i)] $x$ is $\zeta$-regular for $\gamma G$ if and only if $\alpha x$ is $\alpha \zeta$-regular for $\alpha \gamma G$
\item[(ii)] $\{\zeta \mid \zeta$ is regular for $\alpha \gamma
                G\}=\{\alpha \zeta \mid \zeta$ is regular for $\gamma G\}.$
                \end{itemize}
                \end{lemma}
                \begin{proof}[{\em Proof of Lemma \ref{lemma:cosetalpha}}]
\begin{itemize}
\item[(i)] The element $x \in \gamma G$ is $\zeta$-regular
          \begin{align*}
          &\Leftrightarrow V
                  (x, \zeta)\mbox{ contains a regular vector $v$}\\
          &\Leftrightarrow V
                  (\alpha x, \alpha \zeta)\mbox{ contains a regular vector $v$}\\
          &\Leftrightarrow
                  \alpha x \in \alpha \gamma G\mbox{ is $\alpha \zeta$-regular
                  for $\alpha \gamma G$}.\end{align*}
\item[(ii)] The eigenvalue $\zeta \in \mathbb{C}$
                  is regular for $\gamma G$
\begin{align*}
&\Leftrightarrow V (x,                  \zeta)\mbox{ contains a regular vector $v$ for some $x \in \gamma G$}\\
                  &\Leftrightarrow V (\alpha x, \alpha \zeta)\mbox {
                  contains a regular vector $v$ for some $\alpha x \in \alpha
                  \gamma G$}\\
                  &\Leftrightarrow \alpha \zeta \in \mathbb{C}\mbox{
                  is regular for $\alpha \gamma G$}.
                  \end{align*} 
\end{itemize}
\end{proof}
              
              \begin{corollary}\label{corollary:regularcoset} 
              If $x \in \gamma G$ is
                $\zeta$-regular, then $\alpha x \in \alpha \gamma G$ is
                $\alpha \zeta$-regular (in $\alpha \gamma G)$, and the
                eigenvalues of $\alpha x$ acting on $V$ are $\{\alpha
                \varepsilon_i^{\ast} \zeta^{m_i^{\ast}} \mid i = 1, \ldots,
                n\}.$ (See {\S}\ref{subsection:invariant} for a definition
                of the coexponents $m_i^{\ast}$, and
                {\S}\ref{subsection:coset} for a definition of the
                $V^{\ast}$-factors $\varepsilon_i^{\ast}$.)

              \end{corollary}
              
              \begin{proof}[{\em Proof of Corollary \ref{corollary:regularcoset}}]
                This follows directly from Corollary
                \ref{corollary:coseteigenvalues} and Lemma
                \ref{lemma:cosetalpha}.
              \end{proof}
              
              By Lemma \ref{lemma:eigenvalue} and Corollary
              \ref{corollary:regularcoset}, if we can find $\zeta \in
              \mathbb{C}$ which is regular and satisfies $\alpha
              \varepsilon_i^{\ast} \zeta^{m_i^{\ast}} \neq 1$ for all $i = 1, \ldots, r)$, it will follow that $\tmop{Fix}_V
              (\alpha x) =\{0\}$, where $x \in \gamma G$ is any
              $\zeta$-regular element.

As the proofs for the remaining six cases are very similar, we present only the most difficult case here.  The full proof is contained in \cite{Koonin2012}.

           {\bf Case 2:} $G = G (2, 2, 4)$, $\gamma = \dfrac{1}{2}
            \left(\begin{array}{cccc}
              1 &1 &1 &- 1\\
              1 &1 &1 &1\\
              1 &- 1& 1& 1\\
              1 &- 1 &- 1 &1
            \end{array}\right)$

            $m^{\ast}_i = 1,
            3, 3, 5$; $\varepsilon_i^{\ast} = 1$, $\omega, \omega^2, 1$; $\zeta$ regular for $\gamma G \Leftrightarrow
            \left| \zeta \right| \in \{1, 2, 3, 6, 12\}.$

            We must choose $\zeta$ such that $\alpha \zeta \neq 1$, $\alpha
            \omega \zeta^3 \neq 1$, $\alpha \omega^2 \zeta^3 \neq 1$, and $\alpha \zeta^5 \neq 1$. Upon simplification, this becomes $\zeta
            \neq \alpha^{- 1}$, $\zeta^3 \neq \alpha^{- 1} \omega$, $\alpha^{- 1}
            \omega^2$ and $\zeta^5 \neq \alpha^{- 1}$.

            There are 10 complex numbers $\zeta$ whose order $\left| \zeta \right|
            \in \{1, 2, 3, 6, 12\}$ (all 12th roots of unity, except $i$ and $- i$, which have order 4). Of these, three satisfy
            $\zeta^3 = 1$, three satisfy $\zeta^3 = - 1$, two satisfy $\zeta^3
            = i$ and two satisfy $\zeta^3 = - i$. Hence there are at least
            four 12th roots of unity $\zeta$ (in fact 7, if we argue more
            carefully) which satisfy $\zeta^3 \neq \alpha^{- 1} \omega,
            \alpha^{- 1} \omega^2$. Since all 5th powers of 12th roots of
            unity are distinct (as $\tmop{gcd}(5, 12) = 1$), at most one of these
            satisfies $\zeta^5 = \alpha^{- 1}$, and also at most one satisfies
            $\zeta = \alpha^{- 1}$. Hence we can choose such a $\zeta$ as
            required.

           As mentioned above, the other cases are similar but easier. This completes  the proof of
            Proposition \ref{prop:int6} and Theorem \ref{theorem:poset}.
        \end{proof}
        \end{proof}

    Theorem \ref{theorem:poset} describes the principal upper order ideals of
    the poset $\mathcal{S}_{\zeta}^V (\gamma G) .$ The description of
    principal lower order ideals is much simpler.  Recall (Corollary \ref{corollary:steinberg}) that if $U$ is any subset of $V =\mathbb{C}^n$, then $G_U$ is the parabolic subgroup of $G$ which fixes $U$ pointwise.
    
    \begin{lemma}\label{lemma:lower}
      Let $\gamma G$ be a reflection coset in $V
      =\mathbb{C}^n$. Let $E = V (\gamma, \zeta) \in \mathcal{S}_{\zeta}^V
      (\gamma G)$. Then $\widetilde{\mathcal{S}}_{\zeta}^V
      (\gamma G)_{\leqslant E} \cong \widetilde{\mathcal{S}}^V_{\zeta} (\gamma
      G_E)$
      \end{lemma}
      
    \begin{proof}
        Suppose $V (\gamma h, \zeta) \leqslant E.$ Then $\gamma$ and $\gamma h$ both act on $E$ as multiplication by $\zeta$. Then $h =
        \gamma^{- 1} (\gamma h) \in
        G_E .$

        Conversely, if $h \in G_E$ then certainly $V (\gamma h, \zeta)$
        contains $E$. Note that $G_E$ is itself a unitary reflection group
        by Steinberg's fixed point theorem (Corollary \ref{corollary:steinberg}).
      \end{proof}

\section{Homology of $\widetilde{\mathcal{S}}_{\zeta}^V (\gamma
G)$} \label{subsection:homologyS_v(zeta)}

Recall (Definition \ref{definition:Stilde_v(zeta)}) that
$\widetilde{\mathcal{S}}_{\zeta}^V (\gamma G)$ is the poset obtained from
$\mathcal{S}_{\zeta}^V (\gamma G)$ by removing the unique maximal element, as well as
the unique minimal element if it exists. In this section we assume that $G$ is an
irreducible unitary reflection group.

The main theorem in this section is the following:

\begin{theorem}\label{theorem:bouquet}
Suppose $\gamma G$ is a reflection coset in
  $V$, that $G$ acts irreducibly on $V$, and that $\zeta \in
  \mathbb{C}^{\times}$ is a complex root of unity. Then
  $\widetilde{H}_i ( \widetilde{\mathcal{S}}_{\zeta}^V (\gamma G), \mathbb{Z})
  = 0$ for $i \neq l ( \widetilde{\mathcal{S}}_{\zeta}^V (\gamma
  G))$.
  \end{theorem}

\begin{proof}
    Consider the following propositions:

    \begin{proposition}\label{prop:sphericity1}
      Theorem \ref{theorem:bouquet} holds when $l (
      \widetilde{\mathcal{S}}_{\zeta}^V (\gamma G)) \geqslant 2$.
    \end{proposition}
    
    \begin{proposition}\label{prop:sphericity2}
      Theorem \ref{theorem:bouquet} holds when $l (
      \widetilde{\mathcal{S}}_{\zeta}^V (\gamma G)) \geqslant 2$ and $G$ is
      an exceptional irreducible unitary reflection
      group.
\end{proposition}
    
    Clearly (\ref{theorem:bouquet}) $\Rightarrow$ (\ref{prop:sphericity1})
    $\Rightarrow$ (\ref{prop:sphericity2}). The reverse implications are also
    true:

    \begin{theorem}\label{theorem:sphericity}
    We have the implications (\ref{prop:sphericity2}) $\Rightarrow$ (\ref{prop:sphericity1}) $\Rightarrow$ (\ref{theorem:bouquet}).
    \end{theorem}
    
    \begin{proof}[{\em Proof of Theorem \ref{theorem:sphericity}}]
      (\ref{prop:sphericity1}) $\Rightarrow$ (\ref{theorem:bouquet}). It is
      sufficient to show that Theorem \ref{theorem:bouquet} holds in the case
      $l ( \widetilde{\mathcal{S}}_{\zeta}^V (\gamma G)) \leqslant
      1$.

      The poset $\widetilde{\mathcal{S}}_{\zeta}^V (\gamma G)$ can have
      nonzero homology in dimension at most $l (
      \widetilde{\mathcal{S}}_{\zeta}^V ( \gamma G))$. Hence
      Theorem \ref{theorem:bouquet} is trivially true if $l (
      \widetilde{\mathcal{S}}_{\zeta}^V (\gamma G)) = 0$.

      For the case $l ( \widetilde{\mathcal{S}}_{\zeta}^V (\gamma
      G)) = 1$ we require the following proposition:

      \begin{proposition}\label{prop:connected}
        If $l (\widetilde{\mathcal{S}}_{\zeta}^V (\gamma G)) > 0$, then
        $\widetilde{\mathcal{S}}_{\zeta}^V (\gamma G)$ is
        connected.
        \end{proposition}
        
        \begin{proof}[{\em Proof of Proposition \ref{prop:connected}}]
          Given two eigenspaces $X_1, X_2,$ we must find a sequence of
          eigenspaces $X_1 = Y_1, Y_2, \ldots$, $Y_k = X_2$ such that for each
          $i$, either $Y_i$ covers $Y_{i + 1}$ or $Y_{i + 1}$ covers $Y_i$.
          First note that it suffices to consider the case where $X_1 = E_1$,
          $X_2 = E_2$ are maximal eigenspaces, since every eigenspace is
          contained within a maximal eigenspace.

          We now require a lemma:
          
          \begin{lemma}\label{lemma:eigintersection}
            Let $X = V (\gamma g,
            \zeta)$ be a $\zeta$-eigenspace in $\widetilde{\mathcal{S}}_{\zeta}^V (\gamma G)$,
            \and let $r$ be a reflection in $G$ such that $r \nin N
            (X)$. Then $\dim (r X \cap X) = \dim X - 1$.
            \end{lemma}
            
            \begin{proof}[{\em Proof of Lemma \ref{lemma:eigintersection}}]
              It is clear that $r X \cap X \supseteq \tmop{Fix} (r) \cap X.$ \
              Since $r$ is a reflection, $\tmop{Fix} (r)$ has
              codimension $1$ in $V$. Hence Lemma \ref{lemma:eigintersection} follows.
            \end{proof}
            
            To complete the proof of Proposition \ref{prop:connected}, we know
            (by Theorem \ref{theorem:eigenspacescoset}(ii)) that $E_2 = x
            E_1$ for some $x \in G$. Write $x = r_1 r_2 \ldots r_m$ for some
            reflections $r_i \in G$. If for some $i$, $r_i$ normalises
            $r_{i + 1} r_{i + 2} \ldots r_m E_1$, delete it from the
            expression for $x$. In other words, setting $\hat{x} = r_1 r_2
            \ldots r_{i - 1} r_{i + 1} \ldots r_m$, it remains true that $E_2
            = \hat{x} E_1$. Hence we can assume that in the expression for
            $x$, $r_i$ does not normalise $r_{i + 1} r_{i + 2} \ldots r_m E_1$
            for any $i$. Let $Y_1 = E_1$, $Y_{2 i + 1} = r_{m - i + 1}
            r_{m - i + 2} \ldots r_m E_1$  $(1 \leqslant i \leqslant m)$, and $Y_{2
            i} = Y_{2 i - 1} \cap Y_{2 i + 1}$  $(1 \leqslant i \leqslant m)$. Denote the product $r_{m - i +
            1} r_{m - i + 2} \ldots r_m$ by $x_i$  $(1 \leqslant i \leqslant m)$.
            Thus $x = x_m$. \

            It is clear that, for $1 \leqslant i \leqslant m$, $Y_{2 i + 1} =
            x_i E_1 = x_i V (\gamma g, \zeta) = V (x_i \gamma gx_i^{- 1},
            \zeta) \in \widetilde{\mathcal{S}}_{\zeta}^V (\gamma
            G)$, since $\gamma \in N_{GL(V)} (G)$. Also $\dim
            (Y_{2 i + 1}) = \dim (E_1)$, since $x_i$ is invertible.  Hence
            for $1 \leqslant i \leqslant m, Y_{2 i + 1}$ is a maximal
            eigenspace. Also, by Corollary
            \ref{corollary:joinsemilattice}, $Y_{2 i} = Y_{2 i - 1} \cap
            Y_{2 i + 1} \in \widetilde{\mathcal{S}}_{\zeta}^V (\gamma G)$.
            Note the requirement in the statement of Proposition
            \ref{prop:connected} that $l (
            \widetilde{\mathcal{S}}_{\zeta}^V (\gamma G)) > 0$,
            since otherwise $Y_{2 i}$ is the zero space, and does not lie in
            $\widetilde{\mathcal{S}}_{\zeta}^V (\gamma G)$ by definition.
            Given that $l ( \widetilde{\mathcal{S}}_{\zeta}^V
            (\gamma G)) > 0$, Lemma \ref{lemma:eigintersection} guarantees
            that $Y_{2 i} \neq \{0\},$ and that $Y_{2 i}$ covers both $Y_{2 i - 1}$ and $Y_{2 i + 1}$  $(1 \leqslant i \leqslant m)$.

            Proposition \ref{prop:connected} now follows.
        \end{proof}
        
        To complete the proof that
        (\ref{prop:sphericity1}) $\Rightarrow$ (\ref{theorem:bouquet}), consider
        the case $l ( \widetilde{\mathcal{S}}_{\zeta}^V (\gamma
        G)) = 1$. In this case, nonzero homology can exist only in degrees
        0 and 1. By Proposition \ref{prop:connected}, $\widetilde{H}_0 (
        \widetilde{\mathcal{S}}_{\zeta}^V (\gamma G)) = 0$, and hence
        Theorem \ref{theorem:bouquet} holds. Thus
        (\ref{prop:sphericity1}) $\Rightarrow$ (\ref{theorem:bouquet}), as
        claimed.

        (\ref{prop:sphericity2}) $\Rightarrow$ (\ref{prop:sphericity1}). By
        Theorem \ref{theorem:CMimprimitive}, when $G = G (r, p, n)$ the poset
        $\mathcal{S}_{\zeta}^V (\gamma G (r, p, n))$ is
        Cohen-Macaulay over $\mathbb{Z}$. Hence so is $\widetilde{\mathcal{S}}_{\zeta}^V (\gamma
        G (r, p, n))$, by Corollary \ref{corollary:joinCM}. Thus Theorem
        \ref{theorem:bouquet} holds immediately in these cases. This completes
        the proof of Theorem \ref{theorem:sphericity}. 
    \end{proof}
    
    Hence, in order to prove Theorem \ref{theorem:bouquet}, it suffices to
    prove Proposition \ref{prop:sphericity2}. That is, we may assume in the
    statement of Theorem \ref{theorem:bouquet} that $l (
    \widetilde{\mathcal{S}}_{\zeta}^V (\gamma G)) \geqslant 2$ and that $G$
    is an exceptional irreducible reflection group.

    Let $V$ be a complex vector space of dimension $n$,
    $G$ a reflection group acting irreducibly on $V$, and
    $\gamma G$ be a reflection coset in $V$. See
    {\S}\ref{subsection:coset} for a definition of the factors $\varepsilon_i$
    of $(G, \gamma)$. Recall that for $\zeta \in \mathbb{C}^{\times}$ a root of unity,
    we define
    \[
    A (\zeta) =\{i \mid
    \varepsilon_i \zeta^{d_i} = 1\},
    \]
where $d_i$ is the degree of the basic invariant corresponding to
    $\varepsilon_i .$ Note that there are exactly $n$ basic invariants. Let $a (\zeta) = \left| A
    (\zeta) \right|$.

    \begin{proposition}\label{prop:sphericity3}
      With the notation of the above paragraph, it
      suffices to check Proposition \ref{prop:sphericity2} for all $(G,
      \gamma, V)$ such that $G$ acts irreducibly on $V$, $a (\zeta) \geqslant 3$ and $a (\zeta) \neq n$.
    \end{proposition}
    
    \begin{proof}[{\em Proof of Proposition \ref{prop:sphericity3}}]
      If $l ( \widetilde{\mathcal{S}}_{\zeta}^V (\gamma G)) \geqslant 2,$
      then the dimension of a maximal eigenspace is at least 3. By Theorem
      \ref{theorem:eigenspacescoset}(i) the dimension of a maximal eigenspace
      is equal to $a (\zeta)$ and it follows immediately that $a (\zeta)
      \geqslant 3.$.

      If $a (\zeta) = n$, then applying Theorem
      \ref{theorem:eigenspacescoset}(i) again, the full space $V$ is
      the unique maximal $\zeta$-eigenspace. Thus $V = V (\gamma g, \zeta)$ for
      some $g \in G$. But then for any $V (\gamma h, \zeta) \in
      \widetilde{\mathcal{S}}_{\zeta}^V (\gamma G)$, $V (\gamma h, \zeta) = V
      ((\gamma g)^{- 1} \gamma h, 1) = V (g^{- 1} h, 1)$. Hence
      $\widetilde{\mathcal{S}}_{\zeta}^V (\gamma G) =
      \widetilde{\mathcal{S}}_1^V (G)$, which is a geometric lattice, and therefore is Cohen-Macaulay over $\mathbb{Z}$. Hence we can assume that $a (\zeta)
      \neq n$, and the proof of Proposition \ref{prop:sphericity3} is
      complete. 
    \end{proof}
    
    It remains to check all triples $(G, \gamma, V)$ satisfying the conditions
    of Proposition \ref{prop:sphericity3}. First consider the special case
    $\gamma = \tmop{Id}$. In this case $\varepsilon_i = 1$ for all
$i$.     
    
    A list of irreducible reflection groups and their degrees is
    given, for example, in \cite[Table D.3, p.275]{LeTa2009}. An inspection of
    this table provides the following cases for consideration. Here $\left|
    \zeta \right|$ is the order of $\zeta$, not its modulus.
\begin{itemize}
\item[(i)] $G = K_5, \left|\zeta\right| = 3$
\item[(ii)]  $G = K_5, \left|\zeta\right| = 6$
\item[(iii)] $G = E_6, \left|\zeta\right| = 2$
\item[(iv)] $G = E_6, \left|\zeta\right| = 3$
\item[(v)] $G = E_7, \left|\zeta\right| = 3$
\item[(vi)] $G = E_7, \left|\zeta\right| = 6$
\item[(vii)] $G = E_8, \left|\zeta\right| = 3$
\item[(viii)] $G = E_8, \left|\zeta\right| = 4$
\item[(ix)] $G = E_8, \left|\zeta\right| = 6$
\end{itemize}

    We note immediately that for cases (i) \& (ii), primitive 3rd
    roots of unity are precisely the negatives of primitive 6th
    roots of unity. Also $(- 1) \tmop{Id} \in K_5$, so in fact the posets
    $\widetilde{\mathcal{S}}_{\zeta}^V (G)$ are identical in these cases. The
    same argument applies to cases (v) \& (vi), and (vii) \& (ix).

    As described in the Introduction, the reduced homology of these posets was computed using MAGMA (\cite{BoWiCa1997}, \cite{MAGMA}), the GAP package `Simplicial Homology' (\cite{GAP42008}, \cite{GAP4}, \cite{DuHeSaWe2003}, \cite{SimplicialHomology}), and (for posets (vii), (viii) and (ix)), RedHom (\cite{Redhom}, \cite{CAPD}). In particular, MAGMA was used to generate the posets under consideration, while GAP and RedHom actually performed the homology calculations. Here $\omega$ denotes a primitive 3rd root of unity, while $i$ denotes a primitive 4th root of unity. 
    \begin{align*}
    \widetilde{H}_j ( \widetilde{\mathcal{S}}_{\omega}^{\mathbb{C}^5}
    (K_5), \mathbb{Z}) &=\left\{\begin{array}{ll} \mathbb{Z}^{364},&\mbox{if $j=2$}\\
 0,&\mbox{otherwise}\end{array}\right.\\
	\widetilde{H}_j ( \widetilde{\mathcal{S}}_{- 1}^{\mathbb{C}^6}
    (E_6), \mathbb{Z}) &=\left\{\begin{array}{ll} \mathbb{Z}^{475},&\mbox{if $j=3$}\\
	0,&\mbox{otherwise}\end{array}\right.\\
    \widetilde{H}_j ( \widetilde{\mathcal{S}}_{\omega}^{\mathbb{C}^6}
    (E_6), \mathbb{Z}) &=\left\{\begin{array}{ll} \mathbb{Z}^{649},&\mbox{if $j=2$}\\ 
    0,&\mbox{otherwise}\end{array}\right.\\    
   \widetilde{H}_j ( \widetilde{\mathcal{S}}_{\omega}^{\mathbb{C}^7}
    (E_7), \mathbb{Z}) &=\left\{\begin{array}{ll} \mathbb{Z}^{87\, 751},&\mbox{if $j=2$}\\
    0,&\mbox{otherwise}\end{array}\right.\\
    \widetilde{H}_j ( \widetilde{\mathcal{S}}_{\omega}^{\mathbb{C}^8}
    (E_8), \mathbb{Z}) &=\left\{\begin{array}{ll} \mathbb{Z}^{723\, 681}&\mbox{if $j=3$}\\
    0,&\mbox{otherwise}\end{array}\right.\\
\widetilde{H}_j ( \widetilde{\mathcal{S}}_{i}^{\mathbb{C}^8}
    (E_8), \mathbb{Z}) &=\left\{\begin{array}{ll} \mathbb{Z}^{21\, 888\, 721}&\mbox{if $j=3$}\\
    0,&\mbox{otherwise}\end{array}\right.
  \end{align*}  
    
    In each case, reduced homology is zero except in top possible dimension.   
    This completes the proof of
    Proposition \ref{prop:sphericity2} in the case $\gamma = \tmop{Id}$.

    Now consider the case of arbitrary $\gamma$. In \cite[Table D.5, p.278]{LeTa2009}, there is a list of reflection cosets $(G, \gamma_0)$. Recall
    (see the discussion following Lemma \ref{lemma:innerautomorphism}) that if
    $(G, \gamma)$ is any reflection coset where $G$ acts irreducibly,
    then either $\gamma$ induces an inner automorphism, or else $(G, \gamma)
    = (G, \alpha \gamma_0)$, where $(G, \gamma_0)$ appears in this table, and
    $\alpha \in \mathbb{C}^{\times}$. (See also \cite[Proposition 3.13]{BrMaMi1999}).

    If $\gamma$ induces an inner automorphism, then as in Lemma
    \ref{lemma:innerautomorphism} we may assume that $\gamma$ is a scalar
    multiple of the identity. If $\gamma$ is a scalar multiple of the
    identity, then $\widetilde{\mathcal{S}}_{\zeta}^V (\gamma G) =
    \widetilde{\mathcal{S}}_{\gamma^{- 1} \zeta} (G)$. By the result for the
    $\gamma = \tmop{Id}$ case, Proposition \ref{prop:sphericity2} holds.

    Hence we may assume that $(G, \gamma) = (G, \alpha \gamma_0)$, where $(G,
    \gamma_0)$ can be found in \cite[Table D.5, p.278]{LeTa2009}, and $\alpha
    \in \mathbb{C}^{\times}$. Since $\widetilde{\mathcal{S}}_{\zeta}^V (\alpha
    \gamma_0 G) = \widetilde{\mathcal{S}}_{\alpha^{- 1} \zeta}^V (\gamma_0 G)$ as
    in the paragraph above, it suffices to consider the case $\alpha = 1$.

    Applying Proposition \ref{prop:sphericity3}, and inspecting \cite[Table D.5, p.278]{LeTa2009}, there are only two cases to consider:
    
\begin{alignat*}{2}
\mbox{(i) }\; G &= G (2, 2, 4),\quad&\gamma = \frac{1}{2} \left(\begin{array}{cccc}
      1 &1 &1 &- 1\\
      1 &1 &1 &1\\
      1 &- 1 &1 &1\\
      1 &- 1 &- 1& 1
    \end{array}\right)\\
        \mbox{(ii) }\; G &= G_{28},\quad &\gamma = \dfrac{1}{\sqrt{2}}
    \left(\begin{array}{cccc}
      0& 0 &1 &1\\
      0& 0 &1 &- 1\\
      1& 1 &0 &0\\
      1 &- 1 &0 &0
    \end{array}\right).
\end{alignat*}    
    
    The first of these cases is covered by Theorem \ref{theorem:CMimprimitive}, whose
proof appears	 in \cite{Koonin2012-2}. For completeness, the proof is
included here.

    In case (i), we have $d_i = 2, 4, 4, 6$  and  $\varepsilon_i = 1, \omega, \omega^2, 1.$  For a counterexample in this case it would be necessary to find a $\zeta$ such
    that exactly 3 of $\{\zeta^2, \zeta^4\omega, \zeta^4
    \omega^2, \zeta^6 \}$ are equal to 1. Certainly $\zeta^4 \omega$ and
    $\zeta^4 \omega^2$ are always unequal. This forces $\zeta^2 = 1$, in
    which case neither $\zeta^4 \omega$ nor $\zeta^4 \omega^2$ are equal
    to 1. Hence no such counterexample exists in this case.

    For (ii), we have $d_i = 2, 6, 8, 12$  and  $\varepsilon_i = 1, -1, 1, -1.$  For a counterexample here we would need to find a $\zeta$ such that exactly 3 of
    $\{\zeta^2, - \zeta^6, \zeta^8, - \zeta^{12} \}$ are equal to 1. Now if
    $\zeta^2 = 1$ then both $- \zeta^6$ and $- \zeta^{12}$ are equal to -1.
    Hence $- \zeta^6 = 1$ and $\zeta^8 = 1$. Thus $\zeta^2 = - 1.$ But
    then $- \zeta^{12} = - 1$, and so at most two of $\{\zeta^2, - \zeta^6, \zeta^8, - \zeta^{12} \}$ are equal to 1.
    Thus no counterexample can be found here either. This completes the
    proof of Proposition \ref{prop:sphericity2} and Theorem
    \ref{theorem:bouquet}.
  \end{proof}

\section{The CM Property for $\mathcal{S}_{\zeta}^V (\gamma G)$}

In this section we examine the Cohen-Macaulay property for the poset
$\mathcal{S}_{\zeta}^V (\gamma G)$, and prove Theorem \ref{theorem:CMall}.

The first step in the proof of Theorem \ref{theorem:CMall} is a reduction to the case where $V$ is irreducible as a $G$-module. Although Theorem \ref{theorem:CMall} concerns reflection groups, we will make use of the following lemma for reflection cosets. $G, \gamma$ and $V$ all have their usual meaning.

\begin{lemma}\label{lemma:CMproduct}
If $\left.\mathcal{S}\right._\zeta^V (\gamma G)$ is Cohen-Macaulay over
  $\mathbb{Z}$ for all $(G, \gamma, V)$ such that $V$ is irreducible as a $\widehat{G} = \langle G, \gamma \rangle$-module, then $\left.\mathcal{S}\right._\zeta^V (\gamma G)$ is Cohen-Macaulay over $\mathbb{Z}$ for all $(G, \gamma, V)$.
\end{lemma}

\begin{proof}
Suppose $V = V_1 \oplus \ldots \oplus
    V_k$, where each $V_i$ is an irreducible $\langle G, \gamma \rangle$-module. Then $G
    = G_1 \times \cdots \times G_k$, and $\gamma = \gamma_1 \oplus \cdots
    \oplus \gamma_k$, where $G_i$ acts a reflection group in $V_i$ and
    trivially on $V_j$ for $j \neq i$, and $\gamma_i \in GL(V_i)$
    normalises $G_i$. Suppose $g = (g_1, g_2, \ldots, g_k) \in G$, and $v = v_1 + \cdots + v_k \in V$. Denote by $\gamma_i'$ the
    restriction of $\gamma_i$ to $V_i$. Then $v \in V (\gamma g, \zeta)$ if
    and only if $v_i \in V_i (\gamma_i' g_i, \zeta)$ for each $i$. Hence
    $\mathcal{S}_{\zeta}^V (\gamma G) \cong \mathcal{S}_{\zeta}^{V_1}
    (\gamma_1 G_1) \times \cdots \times \mathcal{S}_{\zeta}^{V_k} (\gamma_k
    G_k)$. The result now follows from \cite[Theorem 7.1]{Baclawski1980} (thanks to Alex Miller for pointing out this theorem\footnotemark \footnotetext{Theorem \ref{theorem:CMall} was originally proved over $\mathbb{F}$. The use of this theorem of Baclawski affords a proof over $\mathbb{Z}$.}).
\end{proof}

Noting that when $\gamma = \tmop{Id}$, the proof of Lemma \ref{lemma:CMproduct} never leaves the class of reflection groups, Theorem \ref{theorem:CMall} will follow if we can show that $\mathcal{S}_{\zeta}^V (G)$ is Cohen-Macaulay over $\mathbb{Z}$ whenever $V$ is irreducible as a $G$-module.

By Theorem \ref{theorem:CMimprimitive}, it suffices to show that $\mathcal{S}_{\zeta}^V (G)$ is Cohen-Macaulay over $\mathbb{Z}$ whenever $G$ is an exceptional irreducible reflection group.

By Corollary \ref{corollary:joinCM}, it is equivalent to show that $\widetilde{\mathcal{S}}_{\zeta}^V (G)$ is Cohen-Macaulay over $\mathbb{Z}$ for all such $G$.

To do this, we make use of Lemma \ref{lemma:garst}. The statement that $\widetilde{H}_i ( \widetilde{\mathcal{S}}_{\zeta}^V (G), \mathbb{Z})
  = 0$ for $i \neq l ( \widetilde{\mathcal{S}}_{\zeta}^V (G))$ is just the case $\gamma = \tmop{Id}$ of Theorem \ref{theorem:bouquet}.  Now choose $X = V (h, \zeta) \in \widetilde{\mathcal{S}}_{\zeta}^V
    (G)$ and consider $\tmop{lk} (X) = \widetilde{\mathcal{S}}_{\zeta}^V
    (G)_{< X} \ast \widetilde{\mathcal{S}}_{\zeta}^V (G)_{> X}$. \
    By Lemma \ref{lemma:garst} it will suffice to show that $\tmop{lk} (X)$ is CM. By Proposition \ref{prop:CMjoin} it is enough to
    show that each of $\widetilde{\mathcal{S}}_{\zeta}^V (G)_{< X}
    \tmop{and} \widetilde{\mathcal{S}}_{\zeta}^V (G)_{> X}$ is CM.

Examining $\widetilde{\mathcal{S}}_{\zeta}^V (G)_{> X}$ first, we know that $X$
    is contained in some maximal eigenspace $E$. Hence
    $\mathcal{S}_{\zeta}^V (G)_{\geqslant X}$ is a principal upper
    order ideal of $[E, \hat{1}]$, which is a geometric
    lattice by Corollary \ref{corollary:S_v(zeta)upper}. Hence
    $\mathcal{S}_{\zeta}^V (G)_{\geqslant X}$ is a geometric lattice
    too (see \cite[Lecture 3]{Stanley2007}), and so is certainly CM.
    Removing unique minimal and maximal elements from a poset does not affect whether
    the poset is CM or not, by Corollary \ref{corollary:joinCM}. Thus
    $\widetilde{\mathcal{S}}_{\zeta}^V (G)_{> X}$ is CM too.

As for $\widetilde{\mathcal{S}}_{\zeta}^V (G)_{\leqslant X}$, a similar argument to that used in Lemma \ref{lemma:lower} shows that $\widetilde{\mathcal{S}}_{\zeta}^V
      (G)_{\leqslant X} \cong \widetilde{\mathcal{S}}^V_{\zeta} (h G_X) \cong \widetilde{\mathcal{S}}^{V / X}_{\zeta} (h G_X) $.  The last isomorphism holds because every element of $h G_X$ actis by multiplication by $\zeta$ on $X$. Note that $h G_X$ is a genuine reflection coset, not just a refection group, since $g$ acts by multiplication by $\zeta$ on $X$.

Now if $l(\widetilde{\mathcal{S}}_{\zeta}^V (G)_{\leqslant X}) \leqslant 1$ then $\widetilde{\mathcal{S}}_{\zeta}^V (G)_{< X}$ is trivially Cohen-Macaulay.  

If $l(\widetilde{\mathcal{S}}_{\zeta}^V (G)_{\leqslant X}) = 2$ then $l(\widetilde{\mathcal{S}}_{\zeta}^V (G)_{< X}) = 1$, and  $\widetilde{\mathcal{S}}_{\zeta}^V (G)_{< X}$ is Cohen-Macaulay since it is connected, by Proposition \ref{prop:connected}.

A list of irreducible reflection groups and their degrees is
    given in \cite[Table D.3, p.275]{LeTa2009}. An inspection of
    this table reveals that the only other possibility is $l(\widetilde{\mathcal{S}}_{\zeta}^{V / X} (G)_{\leqslant X}) = 3$, when $X$ is an eigenspace of dimension $1$, and maximal eigenspaces have dimension $4$.  This occurs for only three posets: $\widetilde{\mathcal{S}}_{- 1}^{\mathbb{C}^6} (E_6)$, $\widetilde{\mathcal{S}}_{\omega}^{\mathbb{C}^8} E_8)$ and $\widetilde{\mathcal{S}}_{i}^{\mathbb{C}^8} (E_8)$, where $\omega$ and $i$ are respectively primitive $3$rd and $4$th roots of unity.
 
Now if $V / X$ is reducible as a $\widehat{G} = \langle G, h \rangle$-module, then Lemma \ref{lemma:CMproduct} shows that $\widetilde{\mathcal{S}}^{V / X}_{\zeta} (h G_X)$ is the product of smaller posets of length less than or equal to $2$, which must be Cohen Macaulay by the arguments used above.  Hence $\widetilde{\mathcal{S}}^{V / X}_{\zeta} (h G_X)$ is Cohen-Macaulay, by \cite[Theorem 7.1]{Baclawski1980}. Thus we may assume that $V / X$ is irreducible as a $\widehat{G} = \langle G, h \rangle$-module.

If $V / X$ is irreducible as a $G$-module, then $\widetilde{H}_i (\widetilde{\mathcal{S}}^{V / X}_{\zeta} (h G_X), \mathbb{Z}) = 0$ if $i \neq l(\widetilde{\mathcal{S}}^{V / X}_{\zeta} (h G_X))$, by Theorem \ref{theorem:bouquet}. Lower order ideals of $\widetilde{\mathcal{S}}^{V / X}_{\zeta} (h G_X)$ have length less than or equal to 2, and hence are Cohen-Macaulay by the arguments used above. Upper order ideals of $\widetilde{\mathcal{S}}^{V / X}_{\zeta} (h G_X)$ are also Cohen-Macaulay by Corollary \ref{corollary:S_v(zeta)upper}. Hence if $V / X$ is irreducible as a $G$-module, then $\widetilde{\mathcal{S}}^{V / X}_{\zeta} (h G_X)$ is Cohen-Macaulay as required.

So assume $\widetilde{\mathcal{S}}^{V / X}_{\zeta} (h G_X)$ is reducible as a $G$-module.  For the three cases we are considering, $V / X$ has dimenion either $5$ or $7$. Since $V / X$ is the direct sum of isomorphic vector spaces $V_i$, each of these vector spaces must have dimension $1$.

But then by Lemma \ref{lemma:Ggammairreducible}. maximal eigenspaces in $\widetilde{\mathcal{S}}^{V / X}_{\zeta} (h G_X)$ must have dimension at most $1$, which contradicts the assumption that $l(\widetilde{\mathcal{S}}_{\zeta}^{V / X} (G)_{\leqslant X}) = 3$.

This completes the proof of Theorem \ref{theorem:CMall}.

\begin{remark}
The question as to whether $\mathcal{S}_{\zeta}^V (\gamma G)$ is
  Cohen-Macaulay over $\mathbb{Z}$ for all reflection cosets $\gamma G$ is open. By Lemma \ref{lemma:CMproduct}, it would suffice to prove this result for all $(G, \gamma, V)$ such that $V$ is irreducible as a $\widehat{G} = \langle G, \gamma \rangle$-module, but this is unknown at this stage.
\end{remark}

\section{Acknowledgements}
The author would like to acknowledge his supervisor Gus Lehrer and associate supervisor Anthony Henderson for their encouragement, patience, generosity and enthuiasm.

During preparation of this work, it became evident that Alex Miller had obtained many of the same results independently.  The author would like to acknowledge this
and thank Alex for useful discussions on this topic.

\bibliographystyle{plain}

\bibliography{bibliography2}

\begin{thebibliography}{10}

\bibitem{Arnold1969}
V.I. Arnol'd.
\newblock The cohomology ring of the colored braid group.
\newblock {\em Mat. Zametki}, 5:138--140, 1969.

\bibitem{Baclawski1980}
K.~Baclawski.
\newblock Cohen-{M}acaulay ordered sets.
\newblock {\em J. Algebra}, 63:226--258, 1980.

\bibitem{BlLe2001}
J.~Blair and G.I. Lehrer.
\newblock Cohomology actions and centralisers in unitary reflection groups.
\newblock {\em Proc. London Math. Soc. (3)}, 83:582--604, 2001.

\bibitem{BoLeMi2006}
C.~Bonnaf\'e, G.I. Lehrer, and J.~Michel.
\newblock Twisted invariant theory for reflection groups.
\newblock {\em Nagoya Math. J.}, 182:135--170, 2006.

\bibitem{BoWiCa1997}
W.~Bosma, J.~Cannon, and C.~Playoust.
\newblock The {M}agma algebra system. {I}. {T}he user language.
\newblock {\em J. Symbolic Comput.}, 24(3-4):235--265, 1997.
\newblock Computational algebra and number theory (London, 1993).

\bibitem{Brieskorn1973}
E.~Brieskorn.
\newblock Sur les groups de tresses (d'apr\`es {V}.{I}. {A}rnol'd).
\newblock In {\em Seminaire Bourbaki 1971/2}, volume 317 of {\em Lecture Notes
  in Mathematics}, pages 21--44. Springer, Berlin, 1973.

\bibitem{BrMaMi1999}
M.~Brou\'e, G.~Malle, and J.~Michel.
\newblock Towards {S}petses {I}.
\newblock {\em Transform. Groups}, 4:157--218, 1999.

\bibitem{CAPD}
CAPD.
\newblock \url{http://capd.ii.uj.edu.pl}.

\bibitem{Cohen1976}
A.M. Cohen.
\newblock Finite complex reflection groups.
\newblock {\em Ann. Sci. Ecole Norm. Sup.}, 9(4):379--436, 1976.

\bibitem{DuHeSaWe2003}
J.~G. Dumas, F.~Heckenbach, B.~D. Saunders, and V.~Welker.
\newblock Computing simplicial homology based on efficient smith normal form
  algorithms.
\newblock In {\em Algebra, geometry and software systems}, pages 177--206.
  Springer, 2003.

\bibitem{Folkman1966}
J.~Folkman.
\newblock The homology groups of a lattice.
\newblock {\em J. Math. Mech.}, 15:631--636, 1966.

\bibitem{GAP4}
GAP.
\newblock \url{http://gap-system.org}.

\bibitem{GAP42008}
The GAP~Group.
\newblock {\em {GAP -- Groups, Algorithms, and Programming, Version 4.4.12}},
  2008.

\bibitem{Garst1979}
P.~F. Garst.
\newblock {\em Cohen-Macaulay complexes and group actions}.
\newblock PhD thesis, University of Wisconsin-Madison, 1979.

\bibitem{Koonin2012-4}
J.E. Koonin.
\newblock Homology representations of unitary reflection groups.
\newblock {\em http://arxiv.org/abs/1303.5155}.

\bibitem{Koonin2012-2}
J.E. Koonin.
\newblock Topology of eigenspace posets for imprimitive reflection groups.
\newblock {\em http://arxiv.org/abs/1208.4435}.

\bibitem{Koonin2012}
J.E. Koonin.
\newblock {\em Topology of eigenspace posets for unitary reflection groups}.
\newblock PhD thesis, University of Sydney, 2012.

\bibitem{KoJu2012}
J.E. Koonin and M.~Juda.
\newblock Computational algebraic topology for unitary reflection groups.
\newblock {\em in preparation}.

\bibitem{LeTa2009}
G.~I. Lehrer and D.~E. Taylor.
\newblock {\em Unitary Reflection Groups}.
\newblock Cambridge Univerity Press, 2009.

\bibitem{Lehrer1995}
G.I. Lehrer.
\newblock Poincar\'e polynomials for unitary reflection groups.
\newblock {\em Invent. Math.}, 120:411--425, 1995.

\bibitem{Lehrer2004}
G.I. Lehrer.
\newblock A new proof of {S}teinberg's fixed point theorem.
\newblock {\em Int. Math. Research Notices}, 28:1407--1411, 2004.

\bibitem{LeSp1999}
G.I. Lehrer and T.A. Springer.
\newblock Intersection multiplicities and reflection subquotients of unitary
  reflection groups.
\newblock In {\em Geometric Group Theory Down Under (Canberra 1996)}, pages
  181--193. de Gruyer., Berlin, 1999.

\bibitem{LeSp1999-2}
G.I. Lehrer and T.A. Springer.
\newblock Reflection subquotients of unitary reflection groups.
\newblock {\em Canad. J. Math.}, 51:1175--1193, 1999.

\bibitem{MAGMA}
MAGMA.
\newblock \url{http://magma.maths.usyd.edu.au}.

\bibitem{Malle2006}
G.~Malle.
\newblock Splitting fields for extended complex reflection groups and {H}ecke
  algebras.
\newblock {\em Transform. Groups}, 11(2):195--216, 2006.

\bibitem{OrSo1980}
P.~Orlik and L.~Solomon.
\newblock Combinatorics and topology of complements of hyperplanes.
\newblock {\em Invent. Math.}, 56:167--189, 1980.

\bibitem{OrSo1980-2}
P.~Orlik and L.~Solomon.
\newblock Unitary reflection groups and cohomology.
\newblock {\em Invent. Math}, 59:77--94, 1980.

\bibitem{OrSo1982}
P.~Orlik and L.~Solomon.
\newblock Arrangements defined by unitary reflection groups.
\newblock {\em Math. Ann.}, 261:339--357, 1982.

\bibitem{OrTe1992}
P.~Orlik and H.~Terao.
\newblock {\em Arrangements of Hyperplanes}.
\newblock Springer-Verlag, Berlin, first edition, 1992.

\bibitem{SimplicialHomology}
Simplicial Homology~GAP package.
\newblock \url{http://eecis.udel.edu/~dumas/Homology}.

\bibitem{Quillen1978}
D.~Quillen.
\newblock Homotopy properties of the poset of nontrivial $p$-subgroups of a
  group.
\newblock {\em Adv. in Math.}, 28(2):101--128, 1978.

\bibitem{Redhom}
RedHom.
\newblock \url{http://redhom.ii.uj.edu.pl}.

\bibitem{Rylands1990}
L.~J. Rylands.
\newblock {\em Representations of classical groups on the homology of their
  split buildings}.
\newblock PhD thesis, University of Sydney, 1990.

\bibitem{ShTo1954}
G.C. Shephard and J.A. Todd.
\newblock Finite unitary reflection groups.
\newblock {\em Canad. J. Math.}, 6:274--304, 1954.

\bibitem{Springer1974}
T.~Springer.
\newblock Regular elements of finite reflection groups.
\newblock {\em Invent. Math.}, 25:159--198, 1974.

\bibitem{Stanley1997}
R.P. Stanley.
\newblock {\em Enumerative Combinatorics Volume 1}, volume~62 of {\em Cambridge
  Studies in Advanced Mathematics}.
\newblock Cambridge University Press, 1997.

\bibitem{Stanley2007}
R.P. Stanley.
\newblock An introduction to hyperplane arrangements.
\newblock In {\em Geometric Combinatorics}, volume~13 of {\em IAS/Park City
  Math. Ser.}, pages 389--496. Amer. Math. Soc., Providence, RI, 2007.

\bibitem{Steinberg1964}
R.~Steinberg.
\newblock Differential equations invariant under finite reflection groups.
\newblock {\em Trans. Amer. Math. Soc.}, 112:392--400, 1964.

\bibitem{Wachs2004}
M.L. Wachs.
\newblock Poset topology: tools and applications.
\newblock In {\em Geometric Combinatorics}, volume~13 of {\em IAS/Park City
  Math. Ser.}, pages 497--615. Amer. Math. Soc., Providence, RI, 2007.

\end{thebibliography}

\end{document}